\documentclass[a4paper, 11pt]{amsart}
\usepackage{amssymb}
\usepackage{amsthm}
\usepackage{fullpage}
\usepackage{graphicx}
\usepackage{subfig}
\usepackage{amsmath,amscd}
\usepackage{xypic}
\usepackage{color}
\usepackage{float}
\usepackage{vmargin}
\usepackage[colorlinks,citecolor=blue,filecolor=black,linkcolor=blue,urlcolor=black]{hyperref}
\usepackage[all,cmtip]{xy}

\newcommand{\Z}{\ensuremath{\mathbb Z}}

\newcommand{\R}{\ensuremath{\mathbb R}}

\newcommand{\CP}{\ensuremath{\mathbb {CP}}}

\theoremstyle{plain}
\newtheorem{thm}{Theorem}[section]
\newtheorem*{thm*}{Theorem}

\newtheorem*{cor*}{Corollary}
\newtheorem*{prop*}{Proposition}

\newtheorem*{lemma*}{Lemma}
\newtheorem*{claim*}{Claim}

\theoremstyle{definition}
\newtheorem{defn}[thm]{Definition}

\newtheorem*{exmp*}{Example}
\newtheorem*{defn*}{Definition}
\newtheorem*{rem*}{Remark}
\newtheorem*{note*}{Note}

\begin{document}

\title{A four dimensional hyperbolic link complement in a standard $S^2 \times S^2$} %Title of paper

\author{Hemanth Saratchandran}
%\address{Mathematical Institute, Andrew Wiles Building,
%   Oxford University, OX2 6GG, United Kingdom}
\email{hemanth.saratchandran@maths.ox.ac.uk}
%\thanks{A. Betancourt is supported by the Mexican Council of Science and Technology
%(CONACyT)}
\date{\today}

\maketitle %\maketitle must follow title, authors, abstract and \pacs
\parskip=0.2cm
\parindent=0.0cm

\begin{abstract}
Using techniques from the theory of Kirby calculus we give an explicit construction of a four dimensional hyperbolic link complement in a
4-manifold that is diffeomorphic to a standard $S^2 \times S^2$. 
\end{abstract}

\tableofcontents

\section{Introduction}

This paper continues the study initiated in \cite{sarat} and \cite{sarat_2} on explicitly constructing four dimensional hyperbolic link complements, and identifying
their diffeomorphism type. In \cite{sarat_2}  we used the theory of Kirby calculus to produce an explicit example of a four dimensional hyperbolic link
complement in a 4-manifold that was diffeomorphic to the standard 4-sphere. The general procedure for trying to identify the diffeomorphism
type of such a hyperbolic link complement is easy enough to describe. One starts with a non-compact finite volume hyperbolic 4-manifold $M$ with some
number of cusps, and builds a Kirby diagram for this manifold.
Each cusp cross-section of $M$ has the structure of a closed flat 3-manifold. In the case that such a cusp cross-section is an 
$S^1$-fibre bundle over a flat surface one can glue in the associated disk bundle to the cusp cross-section producing a closed 4-manifold, which we call
a ``filling'' and denote by $M_0$. The
original hyperbolic 4-manifold can then be seen to be a codimension two link complement in the filled in 4-manifold. Provided one can understand how
the gluing in of the disk bundles is done on the level of the Kirby diagram it is possible to construct a Kirby diagram for the filling. In general
this approach to constructing a Kirby diagram for the filling will lead to a very complicated Kirby diagram. This is the stage where Kirby calculus enters
the picture. The Kirby diagram for $M_0$ may have many cancelling pairs of handles, the hope is that by carrying out various handle cancellations and handle slides
the diagram of $M_0$ can be reduced to something much simpler. In cases where $M_0$ is simply connected it is possible to reduce the Kirby diagram of $M_0$ 
all the way down to a Kirby diagram of a familiar 4-manifold. 

In this paper we carry out this procedure to obtain a Kirby diagram of a filling that
can be reduced all the way down to the Kirby diagram of a standard $S^2 \times S^2$. The non-compact finite volume hyperbolic 4-manifold we start with
is the one numbered 35 in the Ratcliffe-Tschantz census see \cite{ratcliffe} p.117. It is non-orientable with a total of five cusps. Using work carried out in \cite{sarat} we show how 
to produce a Kirby diagram for the orientable double cover of this manifold. Using the techniques outlined in \cite{sarat_2} we show how to obtain a Kirby diagram for the
filling. From here we employ the exact same Kirby moves we used in \cite{sarat_2} to reduce the diagram of the filling all the way down to one of a standard
$S^2 \times S^2$. The main theorem we prove takes the form:

\begin{thm*}(Theorem \ref{mainthm_2})
There exists a collection $L$ of linked 2-tori embedded in a standard $S^2 \times S^2$ such that the complement $(S^2 \times S^2) - L$
admits a finite volume hyperbolic geometry.
\end{thm*}

This paper relies heavily on the techniques and constructions used in \cite{sarat} and \cite{sarat_2}, and should be seen as a continuation of \cite{sarat_2}.
In view of this we recommend that the reader take a look at
those two papers to get an idea of the type of arguments we will be using. In order to stop this paper becoming overly long we have refrained
from giving explicit details of certain computations. Many of the computations carried out are analogous to the types of computations we carried out
in \cite{sarat_2} and we refer the reader to that paper for the details.

\section*{Acknowledgements}
We would like to thank Andras Juhasz, Marc Lackenby and John Parker for the comments and corrections they gave on an earlier draft of this work.

\section{How to construct the orientable double cover}\label{const}

The main aim of this section is to show the reader how, given a presentation of the fundamental group of a non-orientable Ratcliffe-Tschantz manifold, we can obtain a 
presentation of the orientable double cover.  We will then put this technique to use in the next section for a particular example, showing the reader
how to obtain a Kirby diagram for the orientable double cover.

The principle technique in extracting a presentation for the orientable double cover is an easy case of the \textit{Reidemeister-Schreier} rewriting
process.
Fix a non-orientable Ratcliffe-Tschantz manifold, which from here on in we denote by $M$. We know that the orientable double cover of $M$, denoted 
$\widetilde{M}$, corresponds to an index 2 subgroup of $\pi_1(M)$. In fact if we let
\[ \phi : \pi_1(M) \rightarrow \Z_2 \]
denote the homomorphism that sends the orientation preserving isometries, making up $M$, to $1$ and the orientation reversing isometries to $-1$. Then
we have that $ker(\phi) = \pi_1(\widetilde{M})$, and it is clear that $\pi_1(\widetilde{M})$ is an index 2 subgroup of $\pi_1(M)$.

Given a group $G$ and a finite index subgroup $H$ of $G$ we remind the reader of the definition of a transversal to $H$ in $G$.

\begin{defn}
A transversal to $H$ in $G$ is a subset $T$ of $G$ such that 
\[ G = \bigcup_{t \in T} H\cdot t \]
\end{defn}

In our case it is easy to find a transversal we simply take any orientation reversing isometry, call it $\alpha$. Then it should be clear that
the set $\{1, \alpha \}$ is a transversal to $\pi_1(\widetilde{M})$ in $\pi_1(M)$.
Now, let us fix the finite presentation of $\pi_1(M) = \langle X \vert R \rangle$, recall that the set $X$ is given by the side pairing transformations
that make up $M$. Let $\psi : \pi_1(M) \rightarrow T$ be the map that the sends an element $x \in \pi_1(M)$ to its coset representative. In our case
the definition of $\psi$ can be explicitly given in the form:

\begin{itemize}
\item $\psi(x) = 1$ if $x$ is orientation preserving.

\item $\psi(x) = \alpha$ if $x$ is orientation reversing.
\end{itemize}

Let $\rho : \pi_1(M) \rightarrow \pi_1(\widetilde{M})$ be defined by $\pi(x) = x \cdot \psi(x)^{-1}$. In our case the definition of $\rho$ can also
be explicitly written in the form:

\begin{itemize}
\item $\rho(x) = x$ if $x$ is orientation preserving

\item $\rho(x) = x \cdot \alpha^{-1}$ if $x$ is orientation reversing
\end{itemize}

We then have the following theorem.

\begin{thm}[Reidemeister-Schreier rewriting process]
A presentation of $\pi_1(\widetilde{M})$ is given by $\langle X' \vert R' \rangle$, where 
\[ X' = \{ \rho(tx) \neq 1 \vert t \in T \text{ and } x \in X \} \] 
\[ R' = \{ xrx^{-1} \vert x \in X \text{ and } r \in R \} \]
and the word $xrx^{-1}$ is written in terms of elements of $X'$.
\end{thm}

We would like to remind the reader that the above theorem is actually a very special case of the general rewriting process. As we only need the case
of an index two subgroup (since we are dealing with orientable double covers) we thought it best to restrict to this simple case. For the proof
of the theorem (and its general form) we refer the reader to \cite{johnson} p.106.

%The general form of the rewriting process, for finite index subgroups, is as follows. Let $G = \langle X \vert R \rangle $ and $H$ a finite index subgroup.

\section{Example: Double cover of Manifold 35}

In this section we are going to show the reader how to obtain a Kirby diagram for the orientable double cover of the manifold numbered 35 in the Ratcliffe-Tschantz
census. Each of their manifolds has the hyperbolic 24-cell as a fundamental domain, a self dual 4-dimensional ideal hyperbolic polyhedron. For the basic 
construction of the hyperbolic 24-cell that we will use
we refer the reader to \cite{sarat} sect.2. For more background information on various properties of the 24-cell we recommend the reader consult
\cite{cox} Ch.4 and \cite{kerckhoff} Sect.3.
The reader who is not familiar with the Ratcliffe-Tschantz hyperbolic 4-manifolds is advised to 
look at their paper \cite{ratcliffe}, which is very readable. The reader can also consult \cite{sarat} where we give a very brief introduction to these
manifolds, and explain how they code each of their manifolds. The notation used in this and all other sections is the notation we used in \cite{sarat_2}.

We start by giving the side pairing information for manifold 35, which we will denote by $M$.

The side pairing code for this manifold is \textbf{146928}. The explicit side pairings are given as:
\[\xymatrixcolsep{5pc} \xymatrix{ S_{(+1,+1,0,0)}  \ar[r]^a_{k_{(-1,+1,+1,+1)}} & S_{(-1,+1,0,0)} } \hspace{2cm} \xymatrix{S_{(+1,-1,0,0)}  \ar[r]^b_{k_{(-1,+1,+1,+1)}} & S_{(-1,-1,0,0)} } \]

\[\xymatrixcolsep{5pc} \xymatrix{ S_{(+1,0,+1,0)}  \ar[r]^c_{k_{(+1,+1,-1,+1)}} & S_{(+1,0,-1,0)} } \hspace{2cm} \xymatrix{S_{(-1,0,+1,0)}  \ar[r]^d_{k_{(+1,+1,-1,+1)}} & S_{(-1,0,-1,0)} } \]

\[\xymatrixcolsep{5pc} \xymatrix{ S_{(0,+1,+1,0)}  \ar[r]^e_{k_{(+1,-1,-1,+1)}} & S_{(0,-1,-1,0)} } \hspace{2cm} \xymatrix{S_{(0,+1,-1,0)}  \ar[r]^f_{k_{(+1,-1,-1,+1)}} & S_{(0,-1,+1,0)} } \]

\[\xymatrixcolsep{5pc} \xymatrix{ S_{(+1,0,0,+1)}  \ar[r]^g_{k_{(-1,+1,+1,-1)}} & S_{(-1,0,0,-1)} } \hspace{2cm} \xymatrix{S_{(+1,0,0,-1)}  \ar[r]^h_{k_{(-1,+1,+1,-1)}} & S_{(-1,0,0,+1)} } \]

\[\xymatrixcolsep{5pc} \xymatrix{ S_{(0,+1,0,+1)}  \ar[r]^i_{k_{(+1,-1,+1,+1)}} & S_{(0,-1,0,+1)} } \hspace{2cm} \xymatrix{S_{(0,+1,0,-1)}  \ar[r]^j_{k_{(+1,-1,+1,+1)}} & S_{(0,-1,0,-1)} } \]

\[\xymatrixcolsep{5pc} \xymatrix{ S_{(0,0,+1,+1)}  \ar[r]^k_{k_{(+1,+1,+1,-1)}} & S_{(0,0,+1,-1)} } \hspace{2cm} \xymatrix{S_{(0,0,-1,+1)}  \ar[r]^l_{k_{(+1,+1,+1,-1)}} & S_{(0,0,-1,-1)} }. \]
The labelling of the sides of this manifold are given in the following table: \\

\begin{tabular}{|l|l|l||l|l|l|}
 $A$ &  $S_{(+1,+1,0,0)}$ & $(\frac{1}{\sqrt{2}}, \frac{1}{\sqrt{2}}, 0)$  &  $A'$ & $S_{(-1,+1,0,0)}$ & $(\frac{-1}{\sqrt{2}}, \frac{1}{\sqrt{2}}, 0)$  \\
 $B$ &  $S_{(+1,-1,0,0)}$ & $(\frac{1}{\sqrt{2}}, \frac{-1}{\sqrt{2}}, 0)$ &  $B'$ & $S_{(-1,-1,0,0)}$  & $(\frac{-1}{\sqrt{2}}, \frac{-1}{\sqrt{2}}, 0)$ \\
 $C$ &  $S_{(+1,0,+1,0)}$ &  $(\frac{1}{\sqrt{2}}, 0, \frac{1}{\sqrt{2}})$ &  $C'$ & $S_{(+1,0,-1,0)}$ & $(\frac{1}{\sqrt{2}}, 0, \frac{-1}{\sqrt{2}})$ \\
 $D$ &  $S_{(-1,0,+1,0)}$ &  $(\frac{-1}{\sqrt{2}}, 0, \frac{1}{\sqrt{2}})$ &  $D'$ & $S_{(-1,0,-1,0)}$ & $(\frac{-1}{\sqrt{2}}, 0, \frac{-1}{\sqrt{2}})$ \\ 
 $E$ &  $S_{(0,+1,+1,0)}$ &  $(0, \frac{1}{\sqrt{2}} ,\frac{1}{\sqrt{2}})$ &  $E'$ & $S_{(0,-1,-1,0)}$ &  $(0, \frac{-1}{\sqrt{2}} ,\frac{-1}{\sqrt{2}})$ \\ 
 $F$ &  $S_{(0,+1,-1,0)}$ &   $(0, \frac{1}{\sqrt{2}} ,\frac{-1}{\sqrt{2}})$ &  $F'$ & $S_{(0,-1,+1,0)}$ &  $(0, \frac{-1}{\sqrt{2}} ,\frac{1}{\sqrt{2}})$ \\ 
 $G$ &  $S_{(+1,0,0,+1)}$ &  $(1 + \sqrt{2}, 0, 0)$ &  $G'$ & $S_{(-1,0,0,-1)}$ & $(1 - \sqrt{2}, 0, 0)$ \\ 
 $H$ &  $S_{(+1,0,0,-1)}$ &  $(-1 + \sqrt{2}, 0, 0)$ &  $H'$ & $S_{(-1,0,0,+1)}$ & $(-1 - \sqrt{2}, 0, 0)$ \\
 $I$ &  $S_{(0,+1,0,+1)}$ &  $(0, 1 + \sqrt{2}, 0)$ &  $I'$ & $S_{(0,-1,0,+1)}$ &  $(0, -1 - \sqrt{2}, 0)$ \\
 $J$ &  $S_{(0,+1,0,-1)}$ &   $(0, -1 + \sqrt{2}, 0)$ &  $J'$ & $S_{(0,-1,0,-1)}$ &  $(0, 1 - \sqrt{2}, 0)$  \\
 $K$ &  $S_{(0,0,+1,+1)}$ &   $(0, 0, 1 + \sqrt{2})$ &  $K'$ & $S_{(0,0,+1,-1)}$ &  $(0, 0, -1 + \sqrt{2})$   \\
 $L$ &  $S_{(0,0,-1,+1}$ &   $(0, 0, -1 - \sqrt{2})$  &  $L'$ & $S_{(0,0,-1,-1)}$ &  $(0, 0, 1 - \sqrt{2})$  
\end{tabular} \\

The twenty four 2-handles are given in the following table:

\begin{tabular}{|l | l | }
\hline			
1. &  $\xymatrix{
A \cap C \ar[r]^a & A'\cap D\ar[r]^{d} & A'\cap D' \ar[r]^{a^{-1}} & A\cap C' \ar[r]^{c^{-1}} & A\cap C }$  \\ \hline
  
2. &  $\xymatrix{
A \cap E \ar[r]^a & A'\cap E \ar[r]^{e} & B'\cap E' \ar[r]^{b^{-1}} & B\cap E' \ar[r]^{e^{-1}} & A\cap E }$  \\ \hline  
  
3. &  $\xymatrix{
A \cap F \ar[r]^a & A'\cap F\ar[r]^{f} & B'\cap F' \ar[r]^{b^{-1}} & B\cap F' \ar[r]^{f^{-1}} & A\cap F }$  \\ \hline

4. &  $\xymatrix{
A \cap G \ar[r]^a & A'\cap H'\ar[r]^{h^{-1}} & A\cap H \ar[r]^{a} & A'\cap G' \ar[r]^{g^{-1}} & A\cap G }$  \\ \hline

5. &  $\xymatrix{
A \cap I \ar[r]^a & A'\cap I\ar[r]^i & B'\cap I' \ar[r]^{b^{-1}} & B\cap I' \ar[r]^{i^{-1}} & A\cap I }$  \\ \hline

6. &  $\xymatrix{
A \cap J \ar[r]^a & A'\cap J\ar[r]^j & B'\cap J' \ar[r]^{b^{-1}} & B\cap J' \ar[r]^{j^{-1}} & A\cap J }$  \\ \hline  

7. &  $\xymatrix{
B \cap C \ar[r]^b & B'\cap D\ar[r]^d & B'\cap D' \ar[r]^{b^{-1}} & B\cap C' \ar[r]^{c^{-1}} & B\cap C }$  \\ \hline

8. &  $\xymatrix{
B \cap G \ar[r]^b & B'\cap H'\ar[r]^{h^{-1}} & B\cap H \ar[r]^{b} & B'\cap G' \ar[r]^{g^{-1}} & B\cap G }$  \\ \hline

9. &  $\xymatrix{
C \cap E \ar[r]^c & C'\cap F\ar[r]^{f} & C\cap F' \ar[r]^{c} & C'\cap E' \ar[r]^{e^{-1}} & C\cap E }$  \\ \hline

10. &  $\xymatrix{
C \cap G \ar[r]^c & C'\cap G \ar[r]^{g} & D'\cap G' \ar[r]^{d^{-1}} & D\cap G' \ar[r]^{g^{-1}} & C\cap G }$  \\ \hline

11. &  $\xymatrix{
C \cap H \ar[r]^c & C'\cap H \ar[r]^h & D'\cap H' \ar[r]^{d^{-1}} & D\cap H' \ar[r]^{h^{-1}} & C\cap H }$  \\ \hline

12. &  $\xymatrix{
C \cap K \ar[r]^c & C'\cap L \ar[r]^{l} & C'\cap L' \ar[r]^{c^{-1}} & C\cap K' \ar[r]^{k^{-1}} & C\cap K }$  \\   \hline

13 &  $\xymatrix{
D \cap E \ar[r]^d & D'\cap F \ar[r]^{f} & D\cap F' \ar[r]^{d} & D'\cap E' \ar[r]^{e^{-1}} & D\cap E }$  \\   \hline

14. &  $\xymatrix{
D \cap K \ar[r]^d & D'\cap L \ar[r]^{l} & D'\cap L' \ar[r]^{d^{-1}} & D\cap K' \ar[r]^{k^{-1}} & D\cap K }$  \\   \hline

15. &  $\xymatrix{
E \cap I \ar[r]^e & E'\cap I' \ar[r]^{i^{-1}} & F\cap I \ar[r]^{f} & F'\cap I' \ar[r]^{i^{-1}} & E\cap I }$  \\   \hline

16. &  $\xymatrix{
E \cap J \ar[r]^e & E'\cap J' \ar[r]^{j^{-1}} & F\cap J \ar[r]^{f} & F'\cap J' \ar[r]^{j^{-1}} & E\cap J }$  \\   \hline

17. &  $\xymatrix{
E \cap K \ar[r]^e & E'\cap L \ar[r]^{l} & E'\cap L' \ar[r]^{e^{-1}} & E\cap K' \ar[r]^{k^{-1}} & E\cap K }$  \\   \hline

18. &  $\xymatrix{
F \cap L \ar[r]^f & F'\cap K \ar[r]^{k} & F'\cap K' \ar[r]^{f^{-1}} & F\cap L' \ar[r]^{l^{-1}} & F\cap L }$  \\   \hline

19. &  $\xymatrix{
G \cap I \ar[r]^g & G'\cap J \ar[r]^{j} & G'\cap J' \ar[r]^{g^{-1}} & G\cap I' \ar[r]^{i^{-1}} & G\cap I }$  \\   \hline

20. &  $\xymatrix{
G \cap K \ar[r]^g & G'\cap K' \ar[r]^{k^{-1}} & H'\cap K \ar[r]^{h^{-1}} & H\cap K' \ar[r]^{k^{-1}} & G\cap I' }$  \\   \hline

21. &  $\xymatrix{
G \cap L \ar[r]^g & G'\cap L' \ar[r]^{l^{-1}} & H'\cap L \ar[r]^{h^{-1}} & H\cap L' \ar[r]^{l^{-1}} & G\cap L }$  \\   \hline

22. &  $\xymatrix{
H \cap J \ar[r]^h & H'\cap I \ar[r]^{i} & H'\cap I' \ar[r]^{h^{-1}} & H\cap J' \ar[r]^{j^{-1}} & H\cap J }$  \\   \hline

23. &  $\xymatrix{
I \cap K \ar[r]^i & I'\cap K \ar[r]^{k} & J'\cap K' \ar[r]^{j^{-1}} & J\cap K' \ar[r]^{k^{-1}} & I\cap K }$  \\   \hline

24. &  $\xymatrix{
I \cap L \ar[r]^i & I'\cap L \ar[r]^{l} & J'\cap L' \ar[r]^{j^{-1}} & J\cap L' \ar[r]^{l^{-1}} & J\cap L }$  \\   \hline

\end{tabular} \\

Recall that there are ten classes of closed flat 3-manifolds denoted by $\textbf{A}$, $\textbf{B}$, $\textbf{C}$, $\textbf{D}$, $\textbf{E}$,
$\textbf{F}$, $\textbf{G}$, $\textbf{H}$, $\textbf{I}$, and $\textbf{J}$ in the Hantzsche-Wendt notation see \cite{hantzsche}. These are denoted
by $\mathcal{G}_1$, $\mathcal{G}_2$, $\mathcal{G}_3$, $\mathcal{G}_4$, $\mathcal{G}_5$, $\mathcal{G}_6$, $\mathcal{B}_1$, $\mathcal{B}_2$, 
$\mathcal{B}_3$ and $\mathcal{B}_4$ respectively using the notation of Wolf see \cite{wolf} Thm.3.5.5, p.117. The first six are the orientable ones, with
$\textbf{A}$ being the 3-torus, and $\textbf{B}$ being the orientable $S^1$-fibre bundle over the Klein bottle. The last four are all non-orientable.
The manifold $M$ has five cusps with associated boundary components having code $\textbf{GGGGH}$ 
(or in Wolf's notation $\mathcal{B}_1 \mathcal{B}_1 \mathcal{B}_1 \mathcal{B}_1 \mathcal{B}_2$). The non-orientable flat 3-manifolds of type
$\textbf{G}$ and $\textbf{H}$ both have $\textbf{A}$ as their orientable double cover. This means that the orientable double cover of $M$ will also
have five cusps, each one having Euclidean structure type $\textbf{A}$. Furthermore, from the classification theorem it is known that
only $\textbf{A}$, $\textbf{B}$, $\textbf{G}$, $\textbf{H}$, $\textbf{I}$, and $\textbf{J}$ are $S^1$-fibre bundles over a compact flat surface.
In order to carry out a filling we need to work out which prarbolic transformation, in the parabolic subgroup associated to each cusp, corresponds
to the $S^1$-fibre. The way to do this is to use the same methods carried out in \cite{sarat_2} sect.2. 
We will not give the details of the explicit computations of the parabolic subgroups associated to each cusp, the method is exactly
analogous to what we did in \cite{sarat_2} sect.2. Instead, we simply give the following table which outlines the translations we will be filling along.

\begin{tabular}{|l | l | }
\hline	

Ideal vertex & Filling translation \\ \hline

$\{(1,0,0,0), (-1,0,0,0)\}$ & $c$ \\ \hline

$\{(0,1,0,0), (0,-1,0,0)\}$ & $a$ \\ \hline

$\{(0,0,1,0), (0,0,-1,0)\}$ & $k$ \\ \hline

$\{(0,0,0,1), (0,0,0,-1)\}$ & $i$ \\ \hline

 $\{(\pm 1/2, \pm 1/2, \pm 1/2, \pm 1/2)\}$ & $e^{-1}heh^{-1}$ \\ \hline

\end{tabular} \\

Filling along the boundary components corresponding to these five cusps via the above translations algebraically corresponds to adding the translations
as relations to the fundamental group of manifold 35.  It is easy to see that a presentation for the fundamental group of manifold 35 is obtained
by taking generators corresponding to pairs of sides of $P$, and relations given by the 24 codimension 2 equivalence classes. Therefore a presentation for the filling is
obtained by adding the relations $c = a = k = i = e^{-1}heh^{-1} = 1$. One can then simplify the presentation to obtain 
$\langle e, g  \vert e^2, g^2, ege^{-1}g^{-1} \rangle$ (one can do this computation by hand, although it is rather tedious. An easier approach
is to use a computer program, for example using \textbf{Magma} one can input the presentation of the group and then use the \textbf{ReduceGenerators}
command to obtain the simplification), thereby concluding that the filled in manifold has fundamental group $\Z_2 \times \Z_2$. It is then clear that
the orientable double cover has fundamental group $\Z_2$. If we denote the orientable double cover by $\widetilde{M}$, then we have that
its universal cover is a two fold covering, which we will denote by $\widetilde{M}_2$. All the Ratcliffe-Tschantz manifolds have Euler characteristic
1, this implies that $\widetilde{M}_2$ has Euler characteristic 4. Appealing to the classification theorems of \textit{Donaldson} and \textit{Freedman}
we can conclude that the homeomorphism type of $\widetilde{M}_2$ is determined by one of the four manifolds $S^2 \times S^2$, $\CP^2 \# \CP^2$,
$\overline{\CP^2} \# \overline{\CP^2}$ or ${\CP^2} \# \overline{\CP^2}$ . Using the theory of spin structures one can conclude that $\widetilde{M}_2$ must be homeomorphic to
$S^2 \times S^2$. Unfortunately, one cannot conclude anything about the diffeomorphism type of $\widetilde{M}_2$. In fact it is unknown if $S^2 \times S^2$ admits
a unique smooth structure so at this point we cannot rule out the case that we are getting an exotic copy of $S^2 \times S^2$. In order to understand
the diffeomorphism type of $\widetilde{M}_2$ we will resort to the Kirby calculus. We will construct a Kirby diagram
for the orientable double cover $\widetilde{M}$, then apply various elementary moves to simplify the diagram, then take the double cover of this
simplified diagram and conclude that $\widetilde{M}_2$ is in fact diffeomorphic to a standard $S^2 \times S^2$.

The starting point
to obtaining a Kirby diagram for the double cover is to identify the orientation preserving isometries and those that are orientation
reversing. 
Recall, any side pairing transformation is written as the composition $rk$, where $r$ is reflection in the image side and $k$ is a diagonal matrix.
It is clear that $r$ is orientation reversing (since it is a reflection in a hyperplane), therefore we can conclude that the orientation
preserving isometries are those whose $k$-part has an odd number of $-1$'s, and the orientation reversing isometries are those that
have an even number of $-1$'s.
We then find that the orientation preserving isometries are given by $a,b,c,d,i,k$, and the orientation reversing isometries are given by
$e,f,g,h$.

From section \ref{const}, we know that a fundamental domain for the orientable double cover consists of two copies of the 24-cell $P$ 
attached along a codimension one side corresponding to one of the orientation reversing isometries. If we take the orientation reversing isometry
$g^{-1}$ then we can think of the fundamental domain of the orientable double cover to consist of the union $P \cup g^{-1}\cdot P$. In order
to understand how we obtain the double cover from $P \cup g^{-1}\cdot P$ we need to work out what the side pairing transformations are. I.e. which
side gets paired to which side. Before we do this let us comment on a thought that has possibly crossed the readers mind. Why did we pick $g^{-1}$
as opposed to $g$? The reason for this will become apparent soon, but for now let us just say that it turns out to make ones life much easier
when drawing Kirby diagrams if one chooses $g^{-1}$ over $g$.
Appealing to the theory outlined in section \ref{const} we find that if we fix a side pairing transformation $\phi : S \rightarrow S'$, we have two
cases to consider. First, if $\phi$ is orientation preserving we find that:
\[\xymatrixcolsep{4pc} \xymatrix{ S \ar[r]^{\phi} & S' } \text{ and }  \xymatrix{ g^{-1}S \ar[r]^{g^{-1}\phi g} & g^{-1}S' }  \]
and if $\phi$ is orientation reversing then:
\[\xymatrixcolsep{4pc}  \xymatrix{ S \ar[r]^{g^{-1}\phi} & g^{-1}S' } \text{ and }  \xymatrix{ g^{-1}S \ar[r]^{\phi g} & S' }. \]
From here it is easy to work out what the side pairing transformation for the orientable double cover are.
\[\xymatrixcolsep{4pc} \xymatrix{ A  \ar[r]^a & A' } \hspace{0.5cm} \xymatrix{g^{-1}A  \ar[r]^{g^{-1}ag} & g^{-1}A' } \hspace{2cm} \xymatrix{ B  \ar[r]^b & B' }  \hspace{0.5cm} 
 \xymatrix{g^{-1}B  \ar[r]^{g^{-1}bg} & g^{-1}B' } \]

\[\xymatrixcolsep{4pc} \xymatrix{ C  \ar[r]^c & C' } \hspace{0.5cm} \xymatrix{g^{-1}C  \ar[r]^{g^{-1}cg} & g^{-1}C' } \hspace{2cm} \xymatrix{ D  \ar[r]^d & D' }  \hspace{0.5cm} 
 \xymatrix{g^{-1}D  \ar[r]^{g^{-1}dg} & g^{-1}D' } \]

\[\xymatrixcolsep{4pc} \xymatrix{ E  \ar[r]^{g^{-1}e} & g^{-1}E' } \hspace{0.5cm} \xymatrix{g^{-1}E  \ar[r]^{eg} & E' } \hspace{2cm} \xymatrix{ F  \ar[r]^{g^{-1}f} & g^{-1}F' }  \hspace{0.5cm} 
 \xymatrix{g^{-1}F  \ar[r]^{fg} & F' } \]

\[\xymatrixcolsep{4pc} \xymatrix{ G  \ar[r]^{g^{-1}g} & g^{-1}G' } \hspace{0.5cm} \xymatrix{g^{-1}G  \ar[r]^{gg} & G' } \hspace{2cm} \xymatrix{ H  \ar[r]^{g^{-1}h} & g^{-1}H' }  \hspace{0.5cm} 
 \xymatrix{g^{-1}H  \ar[r]^{hg} & H' } \]

\[\xymatrixcolsep{4pc} \xymatrix{ I  \ar[r]^i & I' } \hspace{0.5cm} \xymatrix{g^{-1}I  \ar[r]^{g^{-1}ig} & g^{-1}I' } \hspace{2cm} \xymatrix{ J  \ar[r]^j & J' }  \hspace{0.5cm} 
 \xymatrix{g^{-1}J  \ar[r]^{g^{-1}jg} & g^{-1}J' } \]

\[\xymatrixcolsep{4pc} \xymatrix{ K  \ar[r]^k & K' } \hspace{0.5cm} \xymatrix{g^{-1}K  \ar[r]^{g^{-1}kg} & g^{-1}K' } \hspace{2cm} \xymatrix{ L  \ar[r]^l & L' }  \hspace{0.5cm} 
 \xymatrix{g^{-1}L  \ar[r]^{g^{-1}lg} & g^{-1}L' } \]
 
The next step is to work out equivalence classes of codimension 2 sides, this is done in the usual way (see \cite{sarat} p.17 for a detailed explanation of how to do this
for any of the Ratcliffe-Tschantz manifolds).
One simply takes a codimension 2 side and applies 
side pairing transformations until one cycles back to the original codimension 2 side. A codimension 2 side is given by the
intersection of two distinct codimension 1 sides. We already know which codimension 1 sides in $P$ intersect, applying the transformation $g^{-1}$ then tells us
which codimension 1 sides intersect in $g^{-1}P$. In total we obtain forty eight distinct equivalence classes, with each class containing precisely 
four distinct codimension 2 sides. The following table collects together all forty eight equivalence classes.

\begin{table}[H]
  
  \resizebox{0.80\textwidth}{!}{\begin{minipage}{\textwidth}

\begin{tabular}{|l | l | }

\hline			
1. &  $\xymatrixcolsep{5pc} \xymatrix{
A \cap C \ar[r]^a & A'\cap D\ar[r]^{d} & A'\cap D' \ar[r]^{a^{-1}} & A\cap C' \ar[r]^{c^{-1}} & A\cap C }$  \\ \hline
  
2. &  $\xymatrixcolsep{5pc} \xymatrix{
A \cap E \ar[r]^a & A'\cap E \ar[r]^{g^{-1}e} & g^{-1}B'\cap g^{-1}E' \ar[r]^{(g^{-1}bg)^{-1}} & g^{-1}B\cap g^{-1}E' \ar[r]^{( g^{-1}e)^{-1}} & A\cap E }$  \\ \hline  
  
3. &  $\xymatrixcolsep{5pc}  \xymatrix{
A \cap F \ar[r]^a & A'\cap F\ar[r]^{ g^{-1}f} &  g^{-1}B'\cap  g^{-1}F' \ar[r]^{( g^{-1}bg)^{-1}} &  g^{-1}B\cap  g^{-1}F' \ar[r]^{( g^{-1}f)^{-1}} & A\cap F }$  \\ \hline

4. &  $\xymatrixcolsep{5pc}  \xymatrix{
A \cap G \ar[r]^{a} & A'\cap H'\ar[r]^{(hg)^{-1}} &  g^{-1}A\cap g^{-1}H \ar[r]^{g^{-1}ag} &  g^{-1}A'\cap  g^{-1}G' \ar[r]^{( g^{-1}g)^{-1}} & A\cap G }$  \\ \hline

5. &  $\xymatrixcolsep{5pc}  \xymatrix{
A \cap I \ar[r]^a & A'\cap I\ar[r]^i & B'\cap I' \ar[r]^{b^{-1}} & B\cap I' \ar[r]^{i^{-1}} & A\cap I }$  \\ \hline

6. &  $\xymatrixcolsep{5pc}  \xymatrix{
A \cap J \ar[r]^a & A'\cap J\ar[r]^j & B'\cap J' \ar[r]^{b^{-1}} & B\cap J' \ar[r]^{j^{-1}} & A\cap J }$  \\ \hline  

7. &  $\xymatrixcolsep{5pc}  \xymatrix{
B \cap C \ar[r]^b & B'\cap D\ar[r]^d & B'\cap D' \ar[r]^{b^{-1}} & B\cap C' \ar[r]^{c^{-1}} & B\cap C }$  \\ \hline

8. &  $\xymatrixcolsep{5pc}  \xymatrix{
B \cap G \ar[r]^b & B'\cap H'\ar[r]^{(hg)^{-1}} &  g^{-1}B\cap  g^{-1}H \ar[r]^{g^{-1}bg} &  g^{-1}B'\cap  g^{-1}G' \ar[r]^{( g^{-1}g)^{-1}} & B\cap G }$  \\ \hline

9. &  $\xymatrixcolsep{5pc}  \xymatrix{
C \cap E \ar[r]^c & C'\cap F\ar[r]^{ g^{-1}f} &  g^{-1}C\cap  g^{-1}F' \ar[r]^{g^{-1}cg} &  g^{-1}C'\cap g^{-1}E' \ar[r]^{(g^{-1}e)^{-1}} & C\cap E }$  \\ \hline

10. &  $\xymatrixcolsep{5pc}  \xymatrix{
C \cap G \ar[r]^c & C'\cap G \ar[r]^{ g^{-1}g} &  g^{-1}D'\cap  g^{-1}G' \ar[r]^{( g^{-1}dg)^{-1}} &  g^{-1}D\cap  g^{-1}G' \ar[r]^{( g^{-1}g)^{-1}} & C\cap G }$  \\ \hline

11. &  $\xymatrixcolsep{5pc}  \xymatrix{
C \cap H \ar[r]^c & C'\cap H \ar[r]^{ g^{-1}h} &  g^{-1}D'\cap  g^{-1}H' \ar[r]^{( g^{-1}dg)^{-1}} &  g^{-1}D\cap  g^{-1}H' \ar[r]^{( g^{-1}h)^{-1}} & C\cap H }$  \\ \hline

12. &  $\xymatrixcolsep{5pc}  \xymatrix{
C \cap K \ar[r]^c & C'\cap L \ar[r]^{l} & C'\cap L' \ar[r]^{c^{-1}} & C\cap K' \ar[r]^{k^{-1}} & C\cap K }$  \\   \hline

13 &  $\xymatrixcolsep{5pc}  \xymatrix{
D \cap E \ar[r]^d & D'\cap F \ar[r]^{ g^{-1}f} &  g^{-1}D\cap  g^{-1}F' \ar[r]^{ g^{-1}dg} &  g^{-1}D'\cap  g^{-1}E' \ar[r]^{( g^{-1}e)^{-1}} & D\cap E }$  \\   \hline

14. &  $\xymatrixcolsep{5pc}  \xymatrix{
D \cap K \ar[r]^d & D'\cap L \ar[r]^{l} & D'\cap L' \ar[r]^{d^{-1}} & D\cap K' \ar[r]^{k^{-1}} & D\cap K }$  \\   \hline

15. &  $\xymatrixcolsep{5pc}  \xymatrix{
E \cap I \ar[r]^{ g^{-1}e} &  g^{-1}E'\cap  g^{-1}I' \ar[r]^{( g^{-1}ig)^{-1}} &  g^{-1}F\cap  g^{-1}I \ar[r]^{fg} & F'\cap I' \ar[r]^{i^{-1}} & E\cap I }$  \\   \hline

16. &  $\xymatrixcolsep{5pc}  \xymatrix{
E \cap J \ar[r]^{ g^{-1}e} &  g^{-1}E'\cap  g^{-1}J' \ar[r]^{( g^{-1}jg)^{-1}} &  g^{-1}F\cap  g^{-1}J \ar[r]^{fg} & F'\cap J' \ar[r]^{j^{-1}} & E\cap J }$  \\   \hline

17. &  $\xymatrixcolsep{5pc}  \xymatrix{
E \cap K \ar[r]^{ g^{-1}e} &  g^{-1}E'\cap  g^{-1}L \ar[r]^{ g^{-1}lg} &  g^{-1}E'\cap  g^{-1}L' \ar[r]^{( g^{-1}e)^{-1}} & E\cap K' \ar[r]^{k^{-1}} & E\cap K }$  \\   \hline

18. &  $\xymatrixcolsep{5pc}  \xymatrix{
F \cap L \ar[r]^{ g^{-1}f} &  g^{-1}F'\cap  g^{-1}K \ar[r]^{ g^{-1}kg} &  g^{-1}F'\cap  g^{-1}K' \ar[r]^{( g^{-1}f)^{-1}} & F\cap L' \ar[r]^{l^{-1}} & F\cap L }$  \\   \hline

19. &  $\xymatrixcolsep{5pc}  \xymatrix{
G \cap I \ar[r]^{ g^{-1}g} &  g^{-1}G'\cap  g^{-1}J \ar[r]^{ g^{-1}jg} &  g^{-1}G'\cap  g^{-1}J' \ar[r]^{( g^{-1}g)^{-1}} & G\cap I' \ar[r]^{i^{-1}} & G\cap I }$  \\   \hline

20. &  $\xymatrixcolsep{5pc}  \xymatrix{
G \cap K \ar[r]^{ g^{-1}g} &  g^{-1}G'\cap  g^{-1}K' \ar[r]^{( g^{-1}kg)^{-1}} &  g^{-1}H'\cap  g^{-1}K \ar[r]^{( g^{-1}h)^{-1}} & H\cap K' \ar[r]^{k^{-1}} & G\cap K }$  \\   \hline

21. &  $\xymatrixcolsep{5pc}  \xymatrix{
G \cap L \ar[r]^{ g^{-1}g} &  g^{-1}G'\cap  g^{-1}L' \ar[r]^{( g^{-1}lg)^{-1}} &  g^{-1}H'\cap  g^{-1}L \ar[r]^{( g^{-1}h)^{-1}} & H\cap L' \ar[r]^{l^{-1}} & G\cap L }$  \\   \hline

22. &  $\xymatrixcolsep{5pc}  \xymatrix{
H \cap J \ar[r]^{ g^{-1}h} &  g^{-1}H'\cap  g^{-1}I \ar[r]^{ g^{-1}ig} &  g^{-1}H'\cap  g^{-1}I' \ar[r]^{( g^{-1}h)^{-1}} & H\cap J' \ar[r]^{j^{-1}} & H\cap J }$  \\   \hline

23. &  $\xymatrixcolsep{5pc}  \xymatrix{
I \cap K \ar[r]^i & I'\cap K \ar[r]^{k} & J'\cap K' \ar[r]^{j^{-1}} & J\cap K' \ar[r]^{k^{-1}} & I\cap K }$  \\   \hline

24. &  $\xymatrixcolsep{5pc}  \xymatrix{
I \cap L \ar[r]^i & I'\cap L \ar[r]^{l} & J'\cap L' \ar[r]^{j^{-1}} & J\cap L' \ar[r]^{l^{-1}} & I\cap L }$  \\   \hline

\end{tabular} 
     \end{minipage}}
\end{table}

\begin{table}[H]
  
  \resizebox{0.80\textwidth}{!}{\begin{minipage}{\textwidth}

\hskip-4.0cm\begin{tabular}{|l | l | }

\hline			

25. &  $\xymatrixcolsep{5pc} \xymatrix{
 g^{-1}A \cap  g^{-1}C \ar[r]^{ g^{-1}ag} &  g^{-1}A'\cap  g^{-1}D\ar[r]^{ g^{-1}dg} &  g^{-1}A'\cap  g^{-1}D' \ar[r]^{( g^{-1}ag)^{-1}} &  g^{-1}A\cap  g^{-1}C' \ar[r]^{( g^{-1}cg)^{-1}} &  g^{-1}A\cap  g^{-1}C }$  \\ \hline
  
26. &  $\xymatrixcolsep{5pc} \xymatrix{
 g^{-1}A \cap  g^{-1}E \ar[r]^{ g^{-1}ag} &  g^{-1}A'\cap  g^{-1}E \ar[r]^{eg} & B'\cap E' \ar[r]^{b^{-1}} & B\cap E' \ar[r]^{(eg)^{-1}} &  g^{-1}A\cap  g^{-1}E }$  \\ \hline  
  
27. &  $\xymatrixcolsep{5pc}  \xymatrix{
 g^{-1}A \cap  g^{-1}F \ar[r]^{ g^{-1}ag} &  g^{-1}A'\cap  g^{-1}F\ar[r]^{fg} & B'\cap F' \ar[r]^{b^{-1}} & B\cap F' \ar[r]^{(fg)^{-1}} &  g^{-1}A\cap  g^{-1}F }$  \\ \hline

28. &  $\xymatrixcolsep{5pc}  \xymatrix{
 g^{-1}A \cap  g^{-1}G \ar[r]^{ g^{-1}ag} &  g^{-1}A'\cap  g^{-1}H'\ar[r]^{( g^{-1}h)^{-1}} & A\cap H \ar[r]^{a} & A'\cap G' \ar[r]^{(gg)^{-1}} & gA\cap gG }$  \\ \hline

29. &  $\xymatrixcolsep{5pc}  \xymatrix{
 g^{-1}A \cap  g^{-1}I \ar[r]^{ g^{-1}ag} &  g^{-1}A'\cap  g^{-1}I\ar[r]^{ g^{-1}ig} &  g^{-1}B'\cap  g^{-1}I' \ar[r]^{( g^{-1}bg)^{-1}} &  g^{-1}B\cap  g^{-1}I' \ar[r]^{( g^{-1}ig)^{-1}} &  g^{-1}A\cap  g^{-1}I }$  \\ \hline

30. &  $\xymatrixcolsep{5pc}  \xymatrix{
 g^{-1}A \cap  g^{-1}J \ar[r]^{ g^{-1}ag} &  g^{-1}A'\cap  g^{-1}J\ar[r]^{ g^{-1}jg} &  g^{-1}B'\cap  g^{-1}J' \ar[r]^{( g^{-1}bg)^{-1}} &  g^{-1}B\cap  g^{-1}J' \ar[r]^{( g^{-1}jg)^{-1}} &  g^{-1}A\cap  g^{-1}J }$  \\ \hline  

31. &  $\xymatrixcolsep{5pc}  \xymatrix{
 g^{-1}B \cap  g^{-1}C \ar[r]^{ g^{-1}bg} &  g^{-1}B'\cap  g^{-1}D\ar[r]^{ g^{-1}dg} &  g^{-1}B'\cap  g^{-1}D' \ar[r]^{( g^{-1}bg)^{-1}} &  g^{-1}B\cap  g^{-1}C' \ar[r]^{( g^{-1}cg)^{-1}} &  g^{-1}B\cap  g^{-1}C }$  \\ \hline

32. &  $\xymatrixcolsep{5pc}  \xymatrix{
 g^{-1}B \cap  g^{-1}G \ar[r]^{ g^{-1}bg} &  g^{-1}B'\cap  g^{-1}H'\ar[r]^{( g^{-1}h)^{-1}} & B\cap H \ar[r]^{b} & B'\cap G' \ar[r]^{(gg)^{-1}} &  g^{-1}B\cap  g^{-1}G }$  \\ \hline

33. &  $\xymatrixcolsep{5pc}  \xymatrix{
 g^{-1}C \cap  g^{-1}E \ar[r]^{ g^{-1}cg} &  g^{-1}C'\cap  g^{-1}F\ar[r]^{fg} & C\cap F' \ar[r]^{c} & C'\cap E' \ar[r]^{(eg)^{-1}} &  g^{-1}C\cap  g^{-1}E }$  \\ \hline

34. &  $\xymatrixcolsep{5pc}  \xymatrix{
 g^{-1}C \cap  g^{-1}G \ar[r]^{ g^{-1}cg} &  g^{-1}C'\cap  g^{-1}G \ar[r]^{gg} & D'\cap G' \ar[r]^{d^{-1}} & D\cap G' \ar[r]^{(gg)^{-1}} &  g^{-1}C\cap  g^{-1}G }$  \\ \hline

35. &  $\xymatrixcolsep{5pc}  \xymatrix{
 g^{-1}C \cap  g^{-1}H \ar[r]^{ g^{-1}cg} &  g^{-1}C'\cap  g^{-1}H \ar[r]^{hg} & D'\cap H' \ar[r]^{d^{-1}} & D\cap H' \ar[r]^{(hg)^{-1}} &  g^{-1}C\cap  g^{-1}H }$  \\ \hline

36. &  $\xymatrixcolsep{5pc}  \xymatrix{
 g^{-1}C \cap  g^{-1}K \ar[r]^{ g^{-1}cg} &  g^{-1}C'\cap  g^{-1}L \ar[r]^{ g^{-1}lg} &  g^{-1}C'\cap  g^{-1}L' \ar[r]^{( g^{-1}cg)^{-1}} &  g^{-1}C\cap  g^{-1}K' \ar[r]^{( g^{-1}kg)^{-1}} &  g^{-1}C\cap  g^{-1}K }$  \\   \hline

37. &  $\xymatrixcolsep{5pc}  \xymatrix{
 g^{-1}D \cap  g^{-1}E \ar[r]^{ g^{-1}dg} &  g^{-1}D'\cap  g^{-1}F \ar[r]^{fg} & D\cap F' \ar[r]^{d} & D'\cap E' \ar[r]^{(eg)^{-1}} &  g^{-1}D\cap  g^{-1}E }$  \\   \hline

38. &  $\xymatrixcolsep{5pc}  \xymatrix{
 g^{-1}D \cap  g^{-1}K \ar[r]^{ g^{-1}dg} &  g^{-1}D'\cap  g^{-1}L \ar[r]^{ g^{-1}lg} &  g^{-1}D'\cap  g^{-1}L' \ar[r]^{( g^{-1}dg)^{-1}} &  g^{-1}D\cap  g^{-1}K' \ar[r]^{( g^{-1}kg)^{-1}} &  g^{-1}D\cap  g^{-1}K }$  \\   \hline

39. &  $\xymatrixcolsep{5pc}  \xymatrix{
 g^{-1}E \cap  g^{-1}I \ar[r]^{eg} & E'\cap I' \ar[r]^{i^{-1}} & F\cap I \ar[r]^{ g^{-1}f} &  g^{-1}F'\cap  g^{-1}I' \ar[r]^{( g^{-1}ig)^{-1}} &  g^{-1}E\cap  g^{-1}I }$  \\   \hline

40. &  $\xymatrixcolsep{5pc}  \xymatrix{
 g^{-1}E \cap  g^{-1}J \ar[r]^{eg} & E'\cap J' \ar[r]^{j^{-1}} & F\cap J \ar[r]^{ g^{-1}f} &  g^{-1}F'\cap  g^{-1}J' \ar[r]^{( g^{-1}jg)^{-1}} &  g^{-1}E\cap  g^{-1}J }$  \\   \hline

41. &  $\xymatrixcolsep{5pc}  \xymatrix{
 g^{-1}E \cap  g^{-1}K \ar[r]^{eg} & E'\cap L \ar[r]^{l} & E'\cap L' \ar[r]^{(eg)^{-1}} &  g^{-1}E\cap  g^{-1}K' \ar[r]^{( g^{-1}kg)^{-1}} &  g^{-1}E\cap  g^{-1}K }$  \\   \hline

42. &  $\xymatrixcolsep{5pc}  \xymatrix{
 g^{-1}F \cap  g^{-1}L \ar[r]^{fg} & F'\cap K \ar[r]^{k} & F'\cap K' \ar[r]^{(fg)^{-1}} &  g^{-1}F\cap  g^{-1}L' \ar[r]^{( g^{-1}lg)^{-1}} &  g^{-1}F\cap  g^{-1}L }$  \\   \hline

43. &  $\xymatrixcolsep{5pc}  \xymatrix{
 g^{-1}G \cap  g^{-1}I \ar[r]^{gg} & G'\cap J \ar[r]^{j} & G'\cap J' \ar[r]^{(gg)^{-1}} &  g^{-1}G\cap  g^{-1}I' \ar[r]^{(g^{-1}ig)^{-1}} &  g^{-1}G\cap  g^{-1}I }$  \\   \hline

44. &  $\xymatrixcolsep{5pc}  \xymatrix{
 g^{-1}G \cap  g^{-1}K \ar[r]^{gg} & G'\cap K' \ar[r]^{k^{-1}} & H'\cap K \ar[r]^{hg} &  g^{-1}H\cap  g^{-1}K' \ar[r]^{( g^{-1}kg)^{-1}} &  g^{-1}G\cap  g^{-1}K }$  \\   \hline

45. &  $\xymatrixcolsep{5pc}  \xymatrix{
 g^{-1}G \cap  g^{-1}L \ar[r]^{gg} & G'\cap L' \ar[r]^{l^{-1}} & H'\cap L \ar[r]^{(hg)^{-1}} &  g^{-1}H\cap  g^{-1}L' \ar[r]^{( g^{-1}lg)^{-1}} &  g^{-1}G\cap  g^{-1}L }$  \\   \hline

46. &  $\xymatrixcolsep{5pc}  \xymatrix{
 g^{-1}H \cap  g^{-1}J \ar[r]^{hg} & H'\cap I \ar[r]^{i} & H'\cap I' \ar[r]^{hg} &  g^{-1}H\cap  g^{-1}J' \ar[r]^{( g^{-1}jg)^{-1}} &  g^{-1}H\cap  g^{-1}J }$  \\   \hline

47. &  $\xymatrixcolsep{5pc}  \xymatrix{
 g^{-1}I \cap  g^{-1}K \ar[r]^{ g^{-1}ig} &  g^{-1}I'\cap  g^{-1}K \ar[r]^{ g^{-1}kg} &  g^{-1}J'\cap  g^{-1}K' \ar[r]^{( g^{-1}jg)^{-1}} &  g^{-1}J\cap  g^{-1}K' \ar[r]^{( g^{-1}kg)^{-1}} &  g^{-1}I\cap  g^{-1}K }$  \\   \hline

48. &  $\xymatrixcolsep{5pc}  \xymatrix{
 g^{-1}I \cap  g^{-1}L \ar[r]^{ g^{-1}ig} &  g^{-1}I'\cap  g^{-1}L \ar[r]^{ g^{-1}lg} &  g^{-1}J'\cap  g^{-1}L' \ar[r]^{( g^{-1}jg)^{-1}} &  g^{-1}J\cap  g^{-1}L' \ar[r]^{( g^{-1}lg)^{-1}} &  g^{-1}I\cap  g^{-1}L }$  \\   \hline

\end{tabular} 
     \end{minipage}}
\end{table}  

As was mentioned earlier, a fundamental domain for the orientable double cover $\widetilde{M}$ consists of two copies of $P$ obtained by taking the standard
copy of $P$, and transforming it across the side $G$ via the transformation $g^{-1}$. This gives another copy of $P$ on the other
side of $G$. The union of these two copies of $P$ joined along the side $G$ constitutes a fundamental domain for $\widetilde{M}$. Recall that
the side $G$ corresponds to the sphere with centre $(1,0,0,1)$, and in the handle decomposition the side $G$ corresponds to the point $(1 + \sqrt{2}, 0, 0)$.
Since the transformation $g^{-1} = k_{(-1,1,1,-1)}r$, where $r$ is now reflection in the side $G$,
we can think of that part of the Kirby diagram of $\widetilde{M}$ coming from the $g^{-1}P$ piece as being obtained by taking the part coming from $P$, 
applying the transformation $k_{(-1,1,1,-1)}$ to each component of a 1-handle pair, and then reflecting along a plane parallel to the plane
$y = z = 0$ and lying on the right side of $(1 + \sqrt{2}, 0, 0)$. This reflection along a plane parallel to the plane $y = z = 0$ corresponds
to the $r$ part of the transformation $g^{-1}$. Note that the exact centre of the plane parallel to $y = z = 0$ is not important, for example we can take the
centre to be given by the vector $(3,0,0)$ (as long as it lies to the right of $(1 + \sqrt{2}, 0, 0)$). 
From here it is easy to see how a Kirby diagram for $\widetilde{M}$ will look like. It will consist of the usual Kirby diagram corresponding to
$M$, and then to the right of that part of the one handle labelled $G$ it will consist of a diagram obtained by taking the diagram corresponding to
$M$, applying the transformation $k_{(-1,1,1,-1)}$, and then reflecting through a plane parallel to $y = z = 0$ and centred at $(3,0,0)$. Our procedure
for visualising the Kirby diagram of any one of the Ratcliffe-Tschantz manifolds involved
splitting the total diagram into four diagrams, three such diagrams would correspond to those 2-handles that lie in the x-y, x-z and y-z planes, and
one more diagram corresponded to those 2-handles that did not all lie in such a plane, there were always six such 2-handles (see \cite{sarat_} sect.4). We can similarly 
decompose the Kirby diagram of $\widetilde{M}$ into a collection of four such diagrams. The difference in this case is that each diagram
will have have two components, one coming from that part of the fundamental domain corresponding to $P$, and another coming from that part corresponding
to $g^{-1}P$.

It is time to show the reader how these diagrams look like. When we constructed Kirby diagrams for any one of the Ratcliffe-Tschantz manifolds (see \cite{sarat} p.14), 
each diagram would show six 2-handles with each 2-handle being shown in a particular colour.
The situation now is that each diagram will have twelve 2-handles, hence we will need twelve colours to distinguish each 2-handle. 
So as to avoid confusion right from the start we have included the following table which shows the colours we will be using.

\centerline{\graphicspath{ {Manifold_35/} }\includegraphics[width=7cm, height=8cm]{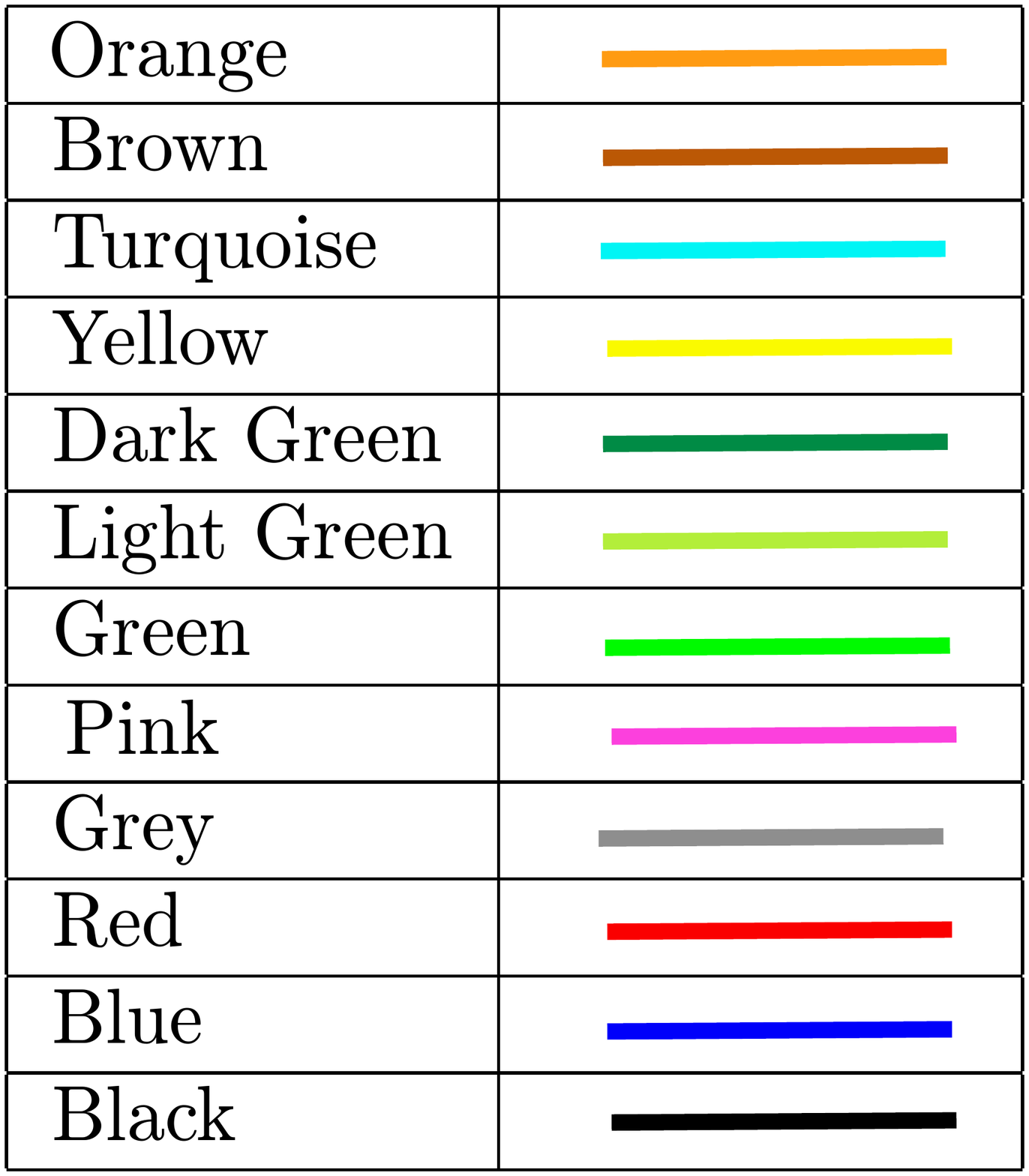}}

In the following diagrams, for each side $S$ we denote the component of the 1-handle corresponding to $g^{-1}S$ by $S-$. There are in total twenty four
1-handles.

The following diagram shows that part of the Kirby diagram contained in the x-y plane, with the table following outlining the colouring of
the 2-handle.

\centerline{\graphicspath{ {Manifold_35/2-handle_cycles/} } \includegraphics[width=12cm, height=6cm]{2cycle_x-y}}

\begin{table}[H]
  
  \resizebox{0.75\textwidth}{!}{\begin{minipage}{\textwidth}

\hskip-4.5cm\begin{tabular}{|l | l | }
\hline
\textbf{Colour} & \textbf{Equivalence class}  \\ \hline

Orange & {$\xymatrixcolsep{5pc}  \xymatrix{
 g^{-1}H \cap  g^{-1}J \ar[r]^{hg} & H'\cap I \ar[r]^{i} & H'\cap I' \ar[r]^{hg} &  g^{-1}H\cap  g^{-1}J' \ar[r]^{( g^{-1}jg)^{-1}} &  g^{-1}H\cap  g^{-1}J }$}  \\ \hline

Brown & {$\xymatrixcolsep{5pc}  \xymatrix{
A \cap G \ar[r]^{a} & A'\cap H'\ar[r]^{(hg)^{-1}} &  g^{-1}A\cap g^{-1}H \ar[r]^{g^{-1}ag} &  g^{-1}A'\cap  g^{-1}G' \ar[r]^{( g^{-1}g)^{-1}} & A\cap G }$}  \\ \hline  

Turquoise & {$\xymatrixcolsep{5pc}  \xymatrix{
 g^{-1}A \cap  g^{-1}J \ar[r]^{ g^{-1}ag} &  g^{-1}A'\cap  g^{-1}J\ar[r]^{ g^{-1}jg} &  g^{-1}B'\cap  g^{-1}J' \ar[r]^{( g^{-1}bg)^{-1}} &  g^{-1}B\cap  g^{-1}J' \ar[r]^{( g^{-1}jg)^{-1}} &  g^{-1}A\cap  g^{-1}J }$}  \\ \hline  

Yellow & {$\xymatrixcolsep{5pc}  \xymatrix{
 g^{-1}A \cap  g^{-1}I \ar[r]^{ g^{-1}ag} &  g^{-1}A'\cap  g^{-1}I\ar[r]^{ g^{-1}ig} &  g^{-1}B'\cap  g^{-1}I' \ar[r]^{( g^{-1}bg)^{-1}} &  g^{-1}B\cap  g^{-1}I' \ar[r]^{( g^{-1}ig)^{-1}} &  g^{-1}A\cap  g^{-1}I }$}  \\ \hline

Dark Green & {$\xymatrixcolsep{5pc}  \xymatrix{
 g^{-1}B \cap  g^{-1}G \ar[r]^{ g^{-1}bg} &  g^{-1}B'\cap  g^{-1}H'\ar[r]^{( g^{-1}h)^{-1}} & B\cap H \ar[r]^{b} & B'\cap G' \ar[r]^{(gg)^{-1}} &  g^{-1}B\cap  g^{-1}G }$}  \\   \hline

Light Green & {$\xymatrixcolsep{5pc}  \xymatrix{
 g^{-1}A \cap  g^{-1}G \ar[r]^{ g^{-1}ag} &  g^{-1}A'\cap  g^{-1}H'\ar[r]^{( g^{-1}h)^{-1}} & A\cap H \ar[r]^{a} & A'\cap G' \ar[r]^{(gg)^{-1}} & gA\cap gG }$ }  \\   \hline

Green & {$\xymatrixcolsep{5pc}  \xymatrix{
B \cap G \ar[r]^b & B'\cap H'\ar[r]^{(hg)^{-1}} &  g^{-1}B\cap  g^{-1}H \ar[r]^{g^{-1}bg} &  g^{-1}B'\cap  g^{-1}G' \ar[r]^{( g^{-1}g)^{-1}} & B\cap G }$ } \\ \hline

Pink & {$\xymatrixcolsep{5pc}  \xymatrix{
G \cap I \ar[r]^{ g^{-1}g} &  g^{-1}G'\cap  g^{-1}J \ar[r]^{ g^{-1}jg} &  g^{-1}G'\cap  g^{-1}J' \ar[r]^{( g^{-1}g)^{-1}} & G\cap I' \ar[r]^{i^{-1}} & G\cap I }$} \\ \hline

Grey & {$\xymatrixcolsep{5pc}  \xymatrix{
 g^{-1}G \cap  g^{-1}I \ar[r]^{gg} & G'\cap J \ar[r]^{j} & G'\cap J' \ar[r]^{(gg)^{-1}} &  g^{-1}G\cap  g^{-1}I' \ar[r]^{(g^{-1}ig)^{-1}} &  g^{-1}G\cap  g^{-1}I }$} \\ \hline

Red & {$\xymatrixcolsep{5pc}  \xymatrix{
A \cap J \ar[r]^a & A'\cap J\ar[r]^j & B'\cap J' \ar[r]^{b^{-1}} & B\cap J' \ar[r]^{j^{-1}} & A\cap J }$} \\ \hline

Blue & {$\xymatrixcolsep{5pc}  \xymatrix{
A \cap I \ar[r]^a & A'\cap I\ar[r]^i & B'\cap I' \ar[r]^{b^{-1}} & B\cap I' \ar[r]^{i^{-1}} & A\cap I }$} \\ \hline

Black & {$\xymatrixcolsep{5pc}  \xymatrix{
H \cap J \ar[r]^{ g^{-1}h} &  g^{-1}H'\cap  g^{-1}I \ar[r]^{ g^{-1}ig} &  g^{-1}H'\cap  g^{-1}I' \ar[r]^{( g^{-1}h)^{-1}} & H\cap J' \ar[r]^{j^{-1}} & H\cap J }$} \\ \hline

\end{tabular} 
      \end{minipage}}
\end{table}

We can see two sets of 1-handles, there are those on the right of $G$ and those on the left. Just to make sure the reader understands exactly
how this diagram is being formed, let us explain why the 1-handle component $A-$ sits where it does. The 1-handle component $A-$ corresponds to the side
$g^{-1}A$, the transformation $g^{-1}$ consists of two parts the r-part, which is reflection through the side $G$ and the $k$-part given by the
diagonal matrix whose diagonal is $(-1, +1, +1, -1)$. The side $A$ can be identified with its centre vector given by $(1,1,0,0)$, which after mapping to $\R^3$
is identified by the co-ordinate $(1/\sqrt{2}, 1/\sqrt{2}, 0)$. Applying
the $k$ matrix to this vector we obtain the vector $(-1,1,0,0)$, this tells us that the $k$ part of the transformation maps the 1-handle component
$A$ to $A'$. We still need to deal with the r-part, since the side $G$ has corresponding 1-handle component in $\R^3$ with centre
$(1+\sqrt{2}, 0,  0)$, and the r-part is reflection through the side $G$, we see that we need to reflect each centre co-ordinate corresponding
to each 1-handle component through a plane parallel to the $y = z = 0$ plane centred at the point $(3,0,0)$ (any centre vector to the right of $G$ will do). From here it 
should be clear that applying $k$ to $A$ followed by reflection in the plane parallel to $y = z = 0$ and centred at $(3,0,0)$ gives us the point where
$A-$ is in the above diagram.  The reader can check for him/her-self that the images of all the other 1-handle components are in the places they are shown
in the diagram.

The fundamental domain for $\widetilde{M}$ consists of two copies of $M$ joined together at the side $G$, therefore the total polyhedron that
constitutes the fundamental domain will not contain the side $G$ and its image $g^{-1}G'$ as these two sides have been identified. This means that
in the Kirby diagram associated to $\widetilde{M}$ we need to kill the 1-handle pair $G-G'-$, we do this by adding a 2-handle whose attaching circle
runs over the 1-handle $G-G'-$ once. It can be seen as the dotted line running from $G$ to $G'-$ in the above diagram, we have chosen to add this 2-handle
so that its attaching circle lies in the x-y and x-z planes.

We move on to show that part of the Kirby diagram that is contained in the x-z plane, the table that follows the diagram shows which 2-handle 
corresponds to which colour.

\centerline{\graphicspath{ {Manifold_35/2-handle_cycles/} } \includegraphics[width=12cm, height=6cm]{2cycle_x-z}}

\begin{table}[H]
  
  \resizebox{0.75\textwidth}{!}{\begin{minipage}{\textwidth}

\hskip-4.5cm\begin{tabular}{|l | l | }
\hline
\textbf{Colour} & \textbf{Equivalence class}  \\ \hline

Orange & {$\xymatrixcolsep{5pc}  \xymatrix{
 g^{-1}C \cap  g^{-1}H \ar[r]^{ g^{-1}cg} &  g^{-1}C'\cap  g^{-1}H \ar[r]^{hg} & D'\cap H' \ar[r]^{d^{-1}} & D\cap H' \ar[r]^{(hg)^{-1}} &  g^{-1}C\cap  g^{-1}H }$}  \\ \hline

Brown & {$\xymatrixcolsep{5pc}  \xymatrix{
C \cap G \ar[r]^c & C'\cap G \ar[r]^{ g^{-1}g} &  g^{-1}D'\cap  g^{-1}G' \ar[r]^{( g^{-1}dg)^{-1}} &  g^{-1}D\cap  g^{-1}G' \ar[r]^{( g^{-1}g)^{-1}} & C\cap G }$}  \\ \hline  

Turquoise & {$\xymatrixcolsep{5pc}  \xymatrix{
 g^{-1}G \cap  g^{-1}K \ar[r]^{gg} & G'\cap K' \ar[r]^{k^{-1}} & H'\cap K \ar[r]^{hg} &  g^{-1}H\cap  g^{-1}K' \ar[r]^{( g^{-1}kg)^{-1}} &  g^{-1}G\cap  g^{-1}K }$}  \\ \hline  

Yellow & {$\xymatrixcolsep{5pc}  \xymatrix{
 g^{-1}D \cap  g^{-1}K \ar[r]^{ g^{-1}dg} &  g^{-1}D'\cap  g^{-1}L \ar[r]^{ g^{-1}lg} &  g^{-1}D'\cap  g^{-1}L' \ar[r]^{( g^{-1}dg)^{-1}} &  g^{-1}D\cap  g^{-1}K' \ar[r]^{( g^{-1}kg)^{-1}} &  g^{-1}D\cap  g^{-1}K }$ }  \\ \hline

Dark Green & {$\xymatrixcolsep{5pc}  \xymatrix{
 g^{-1}G \cap  g^{-1}L \ar[r]^{gg} & G'\cap L' \ar[r]^{l^{-1}} & H'\cap L \ar[r]^{(hg)^{-1}} &  g^{-1}H\cap  g^{-1}L' \ar[r]^{( g^{-1}lg)^{-1}} &  g^{-1}G\cap  g^{-1}L }$}  \\   \hline

Light Green & {$\xymatrixcolsep{5pc}  \xymatrix{
 g^{-1}C \cap  g^{-1}G \ar[r]^{ g^{-1}cg} &  g^{-1}C'\cap  g^{-1}G \ar[r]^{gg} & D'\cap G' \ar[r]^{d^{-1}} & D\cap G' \ar[r]^{(gg)^{-1}} &  g^{-1}C\cap  g^{-1}G }$}  \\   \hline

Green & {$\xymatrixcolsep{5pc}  \xymatrix{
D \cap K \ar[r]^d & D'\cap L \ar[r]^{l} & D'\cap L' \ar[r]^{d^{-1}} & D\cap K' \ar[r]^{k^{-1}} & D\cap K }$} \\ \hline

Pink & {$\xymatrixcolsep{5pc}  \xymatrix{
C \cap K \ar[r]^c & C'\cap L \ar[r]^{l} & C'\cap L' \ar[r]^{c^{-1}} & C\cap K' \ar[r]^{k^{-1}} & C\cap K }$} \\ \hline

Grey & {$\xymatrixcolsep{5pc}  \xymatrix{
 g^{-1}C \cap  g^{-1}K \ar[r]^{ g^{-1}cg} &  g^{-1}C'\cap  g^{-1}L \ar[r]^{ g^{-1}lg} &  g^{-1}C'\cap  g^{-1}L' \ar[r]^{( g^{-1}cg)^{-1}} &  g^{-1}C\cap  g^{-1}K' \ar[r]^{( g^{-1}kg)^{-1}} &  g^{-1}C\cap  g^{-1}K }$} \\ \hline

Red & {$\xymatrixcolsep{5pc}  \xymatrix{
G \cap K \ar[r]^{ g^{-1}g} &  g^{-1}G'\cap  g^{-1}K' \ar[r]^{( g^{-1}kg)^{-1}} &  g^{-1}H'\cap  g^{-1}K \ar[r]^{( g^{-1}h)^{-1}} & H\cap K' \ar[r]^{k^{-1}} & G\cap K }$ } \\ \hline

Blue & {$\xymatrixcolsep{5pc}  \xymatrix{
C \cap H \ar[r]^c & C'\cap H \ar[r]^{ g^{-1}h} &  g^{-1}D'\cap  g^{-1}H' \ar[r]^{( g^{-1}dg)^{-1}} &  g^{-1}D\cap  g^{-1}H' \ar[r]^{( g^{-1}h)^{-1}} & C\cap H }$ } \\ \hline

Black & { $\xymatrixcolsep{5pc}  \xymatrix{
G \cap L \ar[r]^{ g^{-1}g} &  g^{-1}G'\cap  g^{-1}L' \ar[r]^{( g^{-1}lg)^{-1}} &  g^{-1}H'\cap  g^{-1}L \ar[r]^{( g^{-1}h)^{-1}} & H\cap L' \ar[r]^{l^{-1}} & G\cap L }$ } \\ \hline

\end{tabular} 
      \end{minipage}}
\end{table}

The part contained in the y-z plane is shown in the following diagram.

\centerline{\graphicspath{ {Manifold_35/2-handle_cycles/} } \includegraphics[width=12cm, height=6cm]{2cycle_y-z}}

\begin{table}[H]
  
  \resizebox{0.75\textwidth}{!}{\begin{minipage}{\textwidth}

\hskip-4.5cm\begin{tabular}{|l | l | }
\hline
\textbf{Colour} & \textbf{Equivalence class}  \\ \hline

Orange & {$\xymatrixcolsep{5pc}  \xymatrix{
 g^{-1}E \cap  g^{-1}K \ar[r]^{eg} & E'\cap L \ar[r]^{l} & E'\cap L' \ar[r]^{(eg)^{-1}} &  g^{-1}E\cap  g^{-1}K' \ar[r]^{( g^{-1}kg)^{-1}} &  g^{-1}E\cap  g^{-1}K }$}  \\ \hline

Brown & {$\xymatrixcolsep{5pc}  \xymatrix{
F \cap L \ar[r]^{ g^{-1}f} &  g^{-1}F'\cap  g^{-1}K \ar[r]^{ g^{-1}kg} &  g^{-1}F'\cap  g^{-1}K' \ar[r]^{( g^{-1}f)^{-1}} & F\cap L' \ar[r]^{l^{-1}} & F\cap L }$}  \\ \hline  

Turquoise & {$\xymatrixcolsep{5pc}  \xymatrix{
 g^{-1}I \cap  g^{-1}K \ar[r]^{ g^{-1}ig} &  g^{-1}I'\cap  g^{-1}K \ar[r]^{ g^{-1}kg} &  g^{-1}J'\cap  g^{-1}K' \ar[r]^{( g^{-1}jg)^{-1}} &  g^{-1}J\cap  g^{-1}K' \ar[r]^{( g^{-1}kg)^{-1}} &  g^{-1}I\cap  g^{-1}K }$}  \\ \hline  

Yellow & { $\xymatrixcolsep{5pc}  \xymatrix{
 g^{-1}E \cap  g^{-1}J \ar[r]^{eg} & E'\cap J' \ar[r]^{j^{-1}} & F\cap J \ar[r]^{ g^{-1}f} &  g^{-1}F'\cap  g^{-1}J' \ar[r]^{( g^{-1}jg)^{-1}} &  g^{-1}E\cap  g^{-1}J }$}  \\ \hline

Dark Green & {$\xymatrixcolsep{5pc}  \xymatrix{
 g^{-1}E \cap  g^{-1}I \ar[r]^{eg} & E'\cap I' \ar[r]^{i^{-1}} & F\cap I \ar[r]^{ g^{-1}f} &  g^{-1}F'\cap  g^{-1}I' \ar[r]^{( g^{-1}ig)^{-1}} &  g^{-1}E\cap  g^{-1}I }$ }  \\   \hline

Light Green & {$\xymatrixcolsep{5pc}  \xymatrix{
 g^{-1}F \cap  g^{-1}L \ar[r]^{fg} & F'\cap K \ar[r]^{k} & F'\cap K' \ar[r]^{(fg)^{-1}} &  g^{-1}F\cap  g^{-1}L' \ar[r]^{( g^{-1}lg)^{-1}} &  g^{-1}F\cap  g^{-1}L }$}  \\   \hline

Green & {$\xymatrixcolsep{5pc}  \xymatrix{
I \cap L \ar[r]^i & I'\cap L \ar[r]^{l} & J'\cap L' \ar[r]^{j^{-1}} & J\cap L' \ar[r]^{l^{-1}} & I\cap L }$} \\ \hline

Pink & {$\xymatrixcolsep{5pc}  \xymatrix{
E \cap J \ar[r]^{ g^{-1}e} &  g^{-1}E'\cap  g^{-1}J' \ar[r]^{( g^{-1}jg)^{-1}} &  g^{-1}F\cap  g^{-1}J \ar[r]^{fg} & F'\cap J' \ar[r]^{j^{-1}} & E\cap J }$} \\ \hline

Grey & {$\xymatrixcolsep{5pc}  \xymatrix{
 g^{-1}I \cap  g^{-1}L \ar[r]^{ g^{-1}ig} &  g^{-1}I'\cap  g^{-1}L \ar[r]^{ g^{-1}lg} &  g^{-1}J'\cap  g^{-1}L' \ar[r]^{( g^{-1}jg)^{-1}} &  g^{-1}J\cap  g^{-1}L' \ar[r]^{( g^{-1}lg)^{-1}} &  g^{-1}I\cap  g^{-1}L }$ } \\ \hline

Red & {$\xymatrixcolsep{5pc}  \xymatrix{
I \cap K \ar[r]^i & I'\cap K \ar[r]^{k} & J'\cap K' \ar[r]^{j^{-1}} & J\cap K' \ar[r]^{k^{-1}} & I\cap K }$ } \\ \hline

Blue & { $\xymatrixcolsep{5pc}  \xymatrix{
E \cap K \ar[r]^{ g^{-1}e} &  g^{-1}E'\cap  g^{-1}L \ar[r]^{ g^{-1}lg} &  g^{-1}E'\cap  g^{-1}L' \ar[r]^{( g^{-1}e)^{-1}} & E\cap K' \ar[r]^{k^{-1}} & E\cap K }$} \\ \hline

Black & {$\xymatrixcolsep{5pc}  \xymatrix{
E \cap I \ar[r]^{ g^{-1}e} &  g^{-1}E'\cap  g^{-1}I' \ar[r]^{( g^{-1}ig)^{-1}} &  g^{-1}F\cap  g^{-1}I \ar[r]^{fg} & F'\cap I' \ar[r]^{i^{-1}} & E\cap I }$} \\ \hline

\end{tabular} 
      \end{minipage}}
\end{table}

Finally, we have the 2-handles that do not all lie in one of the above planes. There are twelve in total, six coming from each copy of $P$ contributing 
to the fundamental domain (two copies in total).

\centerline{\graphicspath{ {Manifold_35/2-handle_cycles/} } \includegraphics[width=12cm, height=6cm]{2cycle_x-y-z}}

\begin{table}[H]
  
  \resizebox{0.75\textwidth}{!}{\begin{minipage}{\textwidth}

\hskip-4.5cm\begin{tabular}{|l | l | }
\hline
\textbf{Colour} & \textbf{Equivalence class}  \\ \hline

Orange & {$\xymatrixcolsep{5pc}  \xymatrix{
 g^{-1}D \cap  g^{-1}E \ar[r]^{ g^{-1}dg} &  g^{-1}D'\cap  g^{-1}F \ar[r]^{fg} & D\cap F' \ar[r]^{d} & D'\cap E' \ar[r]^{(eg)^{-1}} &  g^{-1}D\cap  g^{-1}E }$ }  \\ \hline

Brown & {$\xymatrixcolsep{5pc} \xymatrix{
A \cap C \ar[r]^a & A'\cap D\ar[r]^{d} & A'\cap D' \ar[r]^{a^{-1}} & A\cap C' \ar[r]^{c^{-1}} & A\cap C }$}  \\ \hline  

Turquoise & {$\xymatrixcolsep{5pc}  \xymatrix{
 g^{-1}C \cap  g^{-1}E \ar[r]^{ g^{-1}cg} &  g^{-1}C'\cap  g^{-1}F\ar[r]^{fg} & C\cap F' \ar[r]^{c} & C'\cap E' \ar[r]^{(eg)^{-1}} &  g^{-1}C\cap  g^{-1}E }$ }  \\ \hline  

Yellow & {$\xymatrixcolsep{5pc} \xymatrix{
 g^{-1}A \cap  g^{-1}C \ar[r]^{ g^{-1}ag} &  g^{-1}A'\cap  g^{-1}D\ar[r]^{ g^{-1}dg} &  g^{-1}A'\cap  g^{-1}D' \ar[r]^{( g^{-1}ag)^{-1}} &  g^{-1}A\cap  g^{-1}C' \ar[r]^{( g^{-1}cg)^{-1}} &  g^{-1}A\cap  g^{-1}C }$}  \\ \hline

Dark Green & {$\xymatrixcolsep{5pc} \xymatrix{
 g^{-1}A \cap  g^{-1}E \ar[r]^{ g^{-1}ag} &  g^{-1}A'\cap  g^{-1}E \ar[r]^{eg} & B'\cap E' \ar[r]^{b^{-1}} & B\cap E' \ar[r]^{(eg)^{-1}} &  g^{-1}A\cap  g^{-1}E }$}  \\   \hline

Light Green & {$\xymatrixcolsep{5pc}  \xymatrix{
 g^{-1}B \cap  g^{-1}C \ar[r]^{ g^{-1}bg} &  g^{-1}B'\cap  g^{-1}D\ar[r]^{ g^{-1}dg} &  g^{-1}B'\cap  g^{-1}D' \ar[r]^{( g^{-1}bg)^{-1}} &  g^{-1}B\cap  g^{-1}C' \ar[r]^{( g^{-1}cg)^{-1}} &  g^{-1}B\cap  g^{-1}C }$}  \\   \hline

Green & { $\xymatrixcolsep{5pc} \xymatrix{
A \cap E \ar[r]^a & A'\cap E \ar[r]^{g^{-1}e} & g^{-1}B'\cap g^{-1}E' \ar[r]^{(g^{-1}bg)^{-1}} & g^{-1}B\cap g^{-1}E' \ar[r]^{( g^{-1}e)^{-1}} & A\cap E }$ } \\ \hline

Pink & {$\xymatrixcolsep{5pc}  \xymatrix{
C \cap E \ar[r]^c & C'\cap F\ar[r]^{ g^{-1}f} &  g^{-1}C\cap  g^{-1}F' \ar[r]^{g^{-1}cg} &  g^{-1}C'\cap g^{-1}E' \ar[r]^{(g^{-1}e)^{-1}} & C\cap E }$} \\ \hline

Grey & {$\xymatrixcolsep{5pc}  \xymatrix{
 g^{-1}A \cap  g^{-1}F \ar[r]^{ g^{-1}ag} &  g^{-1}A'\cap  g^{-1}F\ar[r]^{fg} & B'\cap F' \ar[r]^{b^{-1}} & B\cap F' \ar[r]^{(fg)^{-1}} &  g^{-1}A\cap  g^{-1}F }$ } \\ \hline

Red & {$\xymatrixcolsep{5pc}  \xymatrix{
B \cap C \ar[r]^b & B'\cap D\ar[r]^d & B'\cap D' \ar[r]^{b^{-1}} & B\cap C' \ar[r]^{c^{-1}} & B\cap C }$} \\ \hline

Blue & { $\xymatrixcolsep{5pc}  \xymatrix{
A \cap F \ar[r]^a & A'\cap F\ar[r]^{ g^{-1}f} &  g^{-1}B'\cap  g^{-1}F' \ar[r]^{( g^{-1}bg)^{-1}} &  g^{-1}B\cap  g^{-1}F' \ar[r]^{( g^{-1}f)^{-1}} & A\cap F }$ } \\ \hline

Black & {$\xymatrixcolsep{5pc}  \xymatrix{
D \cap E \ar[r]^d & D'\cap F \ar[r]^{ g^{-1}f} &  g^{-1}D\cap  g^{-1}F' \ar[r]^{ g^{-1}dg} &  g^{-1}D'\cap  g^{-1}E' \ar[r]^{( g^{-1}e)^{-1}} & D\cap E }$} \\ \hline

\end{tabular} 
      \end{minipage}}
\end{table}

Recall that our primary interest is to study, via Kirby calculus, a boundary filling of the manifold $\widetilde{M}$. Before we take this up in the next section
we mention that we have not yet explained how to obtain the 3-handles of the double cover. The procedure is exactly analogous to how we obtained
the 3-handles for the examples considered in \cite{sarat}.
We simply take three distinct codimension 1 sides with non-empty intersection, then apply
side pairing transformations till we cycle back to the original intersection, this constitutes a 3-handle. As we will be dealing with closed
4-manifolds and hence do not have to worry about the 3 and 4-handles (this is due to a theorem of Laudenbach and Po\'enaru, see \cite{gompf} p.116), we will not
be showing tables of the 3-handles nor pictures of how they look like.

\section{Boundary filling of the orientable double cover of Manifold 35.}

In this section we are going to use elementary moves to reduce the Kirby diagram of $\widetilde{M}$, this will help us in identifying
the diffeomorphism type of the double cover of $\widetilde{M}$. in the \cite{sarat_2} sect.4 we explained how all the 2-handles for a filling
had a planar framing, furthermore we explained how the attaching maps being reflections or compositions of reflections with inversion in $S^2$ had the
effect that when we pushed components of 2-handles through attaching spheres of 1-handles nothing ``wild'' could happen i.e. the 2-handle component being
pushed through would not twist around the attaching sphere it came out of. These observations all hold true in the case of the orientable double cover
of manifold 35, and in fact for all of the Ratcliffe-Tschantz manifolds. We will not go through the details of this as they are completely analogous
to what we did for the example considered in \cite{sarat_2} sect.4.

We already mentioned that the boundary type associated to each ideal vertex is given by the code \textbf{GGGGH} (or in \textit{Wolf's} notation 
$\mathcal{B}_1\mathcal{B}_1\mathcal{B}_1\mathcal{B}_1\mathcal{B}_2$), each of these boundaries are themselves non-orientable, therefore in the orientable
double cover they will lift to there own orientable double covers. We also mentioned that we computed a translation in each of the parabolic subgroups
associated to each cusp, the translations we obtained were $c$, $a$, $k$, $i$ and $e^{-1}heh^{-1}$. The first four translations are all given by orientation
preserving transformations, hence in the double cover they correspond to translations in their respective boundary components. Therefore the associated filling
of the corresponding boundary components in the Kirby diagram of $\widetilde{M}$ will consist of adding four 2-handles running over
$C-C'$, $A-A'$, $K-K'$ and $I-I'$ once. The translation 
$e^{-1}heh^{-1}$ can be written as $e^{-1}heh^{-1} = (e^{-1}g)(g^{-1}h)(eg)(g^{-1}h^{-1})$, therefore the corresponding filling of the associated boundary
component in $\widetilde{M}$ will involve adding a 2-handle with four components, 
one running from $E$ to $H$, followed by one running from $H'-$ to $E-$, followed by one running from $E'$ to $H'$, and finally one running from 
$H-$ to $E'-$.

We move on to showing how the Kirby diagrams look with these added 2-handles.
The following diagram shows the x-y plane, by considering a fundamental domain for the ideal vertex class $\{(0,0,0,1), (0,0,0,-1)\}$
we can replace the translation $i$ with $j$, it will be much easier to use the transformation $j$ when we apply elementary moves
to the Kirby diagram, therefore we make this change from now itself. The reader should notice how the added 2-handles running over $A-A'$ and
$J-J'$ lie completely in the x-y plane.

\centerline{\graphicspath{ {Manifold_35/filling_cusps/} }\includegraphics[width=12cm, height=6cm]{2cycle_x-y}}

The following diagram shows the x-z plane, the reader should note how the added 2-handles running over $K-K'$ and $C-C'$ lie entirely in the
x-z plane.

\centerline{\graphicspath{ {Manifold_35/filling_cusps/} }\includegraphics[width=12cm, height=6cm]{2cycle_x-z}}

The following diagram shows a picture of the y-z plane, this plane also contains the added 2-handles running over $J-J'$ and $K-K'$.

\centerline{\graphicspath{ {Manifold_35/filling_cusps/} }\includegraphics[width=12cm, height=6cm]{2cycle_y-z}}

Finally, we have the twelve 2-handles that do not all lie in any one of the above planes.

\centerline{\graphicspath{ {Manifold_35/filling_cusps/} }\includegraphics[width=12cm, height=6cm]{2cycle_x-y-z}}

We have not shown the added 2-handle corresponding to the translation $e^{-1}h^{-1}eh$ that has components running from $E$ to $H$, $H'-$ to $E-$, $E'$ to $H'$ and
$H-$ to $E'-$. Two of the components, the ones running from $E$ to $H$ and $E'$ to $H'$, move from that part of the y-z plane corresponding to the piece of
the Kirby diagram coming from $P$ to the $x-y$ plane, hence they run outside the four diagrams we have been showing. Due to this none of the
elementary moves we carry out to begin with will affect these two components in any way. Similarly, the two components running from $H'-$ to $E-$ and
$H-$ to $E'-$ move from that part of the y-z plane corresponding to the piece of the Kirby diagram coming from $g^{-1}P$ to the piece in the x-y plane, due to
this they will also not be affected by any of the elementary moves to begin with. Therefore we will choose to leave this 2-handle out of our diagrams to start with, 
towards the end when we start doing handle slides that move between planes we will put this 2-handle back in so the reader can see exactly how it is affected.

We are now in the situation where we have various 2-handles that are running over 1-handles once, and hence we have various handle cancelling pairs.
We want to start carrying out several of these cancellations, however we need to be a bit careful when we do so. A few of the 2-handles intersect the other
planes, hence when carrying out cancellations/slides we must keep track of how these intersection points move. We remind the reader of our coding system that
helps keep track of various intersection points. An intersection point in a diagram will be shown via a black dot, the dot will have a code next to it which is
supposed to tell the reader which 2-handle is creating the point of intersection. The code will consist of either four characters or two characters, in the case
that it consists of four characters the first two tell the reader from which plane the 2-handle, creating the intersection point, lies in. The second
two characters tell us which 1-handles the 2-handle runs over, in situations where there are multiple 2-handles running over the 1-handles we will always
make it apparent as to which 2-handle we are talking about. Finally, in the case that the code consists of just two letters we are to immediately take this to mean
that the 2-handle creating the intersection point is residing in the diagram showing the twelve 2-handles that do not all lie in a single plane. The
two characters of the code then tell us which 1-handles this 2-handle is running over.

The following shows pictures of the x-y, x-z and y-z planes respectively, with intersection points added

\centerline{\graphicspath{ {Manifold_35/filling_cusps/} }\includegraphics[width=11cm, height=6cm]{2cycle_x-y_int}}

\centerline{\graphicspath{ {Manifold_35/filling_cusps/} }\includegraphics[width=11cm, height=6cm]{2cycle_x-z_int}}

\centerline{\graphicspath{ {Manifold_35/filling_cusps/} }\includegraphics[width=11cm, height=6cm]{2cycle_y-z_int}}

The first cancellation we are going to carry out is to cancel $G,G'-$ using the black dashed 2-handle in the x-y plane. This will only affect
the diagrams in the x-y and x-y planes.
The following shows how the diagram in the x-y plane changes.

\centerline{\graphicspath{ {Manifold_35/cancelling1_g/} }\includegraphics[width=10cm, height=5cm]{2cycle_x-y_int}}

The diagram in the x-z plane changes as follows.

\centerline{\graphicspath{ {Manifold_35/cancelling1_g/} }\includegraphics[width=11cm, height=5cm]{2cycle_x-z_int}}

The next step we take is to cancel $A,A'$ using the added 2-handle that runs over this 1-handle once. As this 2-handle resides in the
x-y plane and the diagram showing the twelve 2-handles that do not all lie in a single plane, it is only these two diagrams that will be affected.
However, the reader should keep in mind that this 2-handle creates a point of intersection with the y-z plane, hence there will be some changes to the
y-z plane on the level of intersection points.
The following picture shows how the x-y plane changes when we carry out this cancellation.

\centerline{\graphicspath{ {Manifold_35/cancelling2_a/} }\includegraphics[width=11cm, height=5cm]{2cycle_x-y_int}}

Observe that after carrying out the cancellation we have also carried out two handle slides. When we cancel $A,A'$ we get a blue 2-handle component
that loops back into $I$, we push this through $I$ to come out of $I'$, then slide the blue 2-handle off of $I'$ to give a blue 2-handle
running over $B,B'$ once. We also get a red 2-handle component that loops back into $J$, we can perform an analogous slide to obtain a red 2-handle
running over $B,B'$ once. In general, when we carry out such handle cancellations we will also simultaneously carry out handle slides analogous to the one
described above. It should be clear to the reader that we have carried out such handle slides.

The diagram corresponding to the twelve 2-handles that do not all lie in a single plane changes as follows.

\centerline{\graphicspath{ {Manifold_35/cancelling2_a/} }\includegraphics[width=11cm, height=5cm]{2cycle_x-y-z}}

The reader should note that once again we have carried out some handle slides, when we cancel $A,A'$ we obtain a blue 2-handle component that
loops back into $F$, remember that in the double cover our 1-handle pair is $F,F'-$, therefore when we push this blue 2-handle component through $F$ it will come
out of $F'-$, then slide the blue 2-handle off of $F'-$ to obtain a blue 2-handle running over $B-,B'-$ once. Similarly when we cancel $A,A'$ we get a 
green 2-handle component that loops into $E$, and since $E$ is identified to $E'-$ this component can be pushed through to come out of   
$E'-$, we can then slide the green 2-handle into the position shown in the above picture.

The cancellation of $A,A'$ with the added 2-handle that ran over it once causes the intersection point, in the y-z plane, labelled \textbf{XY\_{AA'}}
to disappear, with many new intersection points appearing. The following picture shows the y-z plane with the added intersection points coming from
the above cancellation.

\centerline{\graphicspath{ {Manifold_35/cancelling2_a/} }\includegraphics[width=11cm, height=6cm]{2cycle_y-z_int}}

We move on to cancelling $C,C'$ using the added 2-handle that passes over it once, and that resides in the x-z plane and the diagram corresponding
to the twelve 2-handles that did not all lie in a single plane.

The x-z plane changes as follows:

\centerline{\graphicspath{ {Manifold_35/cancelling3_c/} }\includegraphics[width=11cm, height=5cm]{2cycle_x-z_int}}

The reader should observe that we have also carried out a handle slide. Namely, we have pushed the blue 2-handle component that loops back into $H$
through $H$ and then slid it off of $H'-$ to get a blue 2-handle running over $D-,D'-$ once.

The diagram corresponding to the twelve 2-handles that did not all lie in a single plane changes to the following diagram.

\centerline{\graphicspath{ {Manifold_35/cancelling3_c/} }\includegraphics[width=10cm, height=4.5cm]{2cycle_x-y-z}}

The 2-handle that was used to do the cancellation intersected the x-y plane, hence this point of intersection will disappear with some new ones
coming in place of it. The following picture shows the coding of these new intersection points, it should be clear as to which components of 2-handle
are creating the intersection points.

\centerline{\graphicspath{ {Manifold_35/cancelling3_c/} }\includegraphics[width=10cm, height=5.5cm]{2cycle_x-y_int}}

There were also points of intersection in the y-z plane created by the brown 2-handle component running from $C$ to $D$ and $C'$ to $D'$. 
These will change as follows:

\centerline{\graphicspath{ {Manifold_35/cancelling3_c/} }\includegraphics[width=10cm, height=5cm]{2cycle_y-z_int}}

The next step we take is to cancel $B,B'$ with the red 2-handle in the x-y plane. This cancellation will affect the x-y plane and the diagram
corresponding to the twelve 2-handles that did not all lie in a single plane.

The following picture shows how the x-y plane changes after this cancellation has been carried out. The reader should observe that when we perform this cancellation using 
the red 2-handle, the blue 2-handle component slides to a zero framed unknot hence cancels a 3-handle and can be deleted from the diagram.

\centerline{\graphicspath{ {Manifold_35/cancelling4_b/} }\includegraphics[width=10cm, height=5cm]{2cycle_x-y_int}}

The following picture shows how the diagram corresponding to the twelve 2-handles that did not all lie in a single plane changes.

\centerline{\graphicspath{ {Manifold_35/cancelling4_b/} }\includegraphics[width=10cm, height=5cm]{2cycle_x-y-z}}

The astute reader would have noticed that in the picture above showing the x-y plane, the intersection point labelled \textbf{DD'} has moved
from the right to the left. This is because when we cancelled $B,B'$ in the diagram corresponding to the twelve 2-handles that did not all lie in a single plane 
we also moved the brown 2-handle into the position shown above. As this brown 2-handle was intersecting the x-y plane, this intersection point must also move.

Finally, we note that the red 2-handle used to carry out the above cancellation intersected the y-z plane, hence we will get some new points of intersection
in this plane. Furthermore, the moving of the brown 2-handle we did above will cause the two intersection points in the y-z plane labelled \textbf{DD'} to
disappear.

The following picture shows the structure of the y-z plane after all the above has been carried out.

\centerline{\graphicspath{ {Manifold_35/cancelling4_b/} }\includegraphics[width=10cm, height=5cm]{2cycle_y-z_int}}

We can then cancel $B-,B'-$ with the blue 2-handle in the diagram corresponding to the twelve 2-handles that did not all lie in a single plane. 
The diagram changes as follows.

\centerline{\graphicspath{ {Manifold_35/cancelling5_b-/} }\includegraphics[width=10cm, height=5cm]{2cycle_x-y-z}}

Note that the green 2-handle component that originally ran from $B-$ to $B'-$ slides to a zero framed unknot and can be immediately deleted from
the diagram.

The cancellation also affects the x-y plane, which changes to the following.

\centerline{\graphicspath{ {Manifold_35/cancelling5_b-/} }\includegraphics[width=10cm, height=5cm]{2cycle_x-y_int}}

The blue 2-handle used to carry out the cancellation intersected the y-z plane, it is clear that the intersection points in the y-z plane change
in the following way:

\centerline{\graphicspath{ {Manifold_35/cancelling5_b-/} }\includegraphics[width=10cm, height=5cm]{2cycle_y-z_int}}

So far we have carried out five different handle cancellation moves. The following pictures collect together how the various diagrams
have changed so far. The first picture shows the x-y and x-z planes respectively, and the second shows the y-z plane and the diagram
corresponding to the twelve 2-handles that did not all lie in a single plane.

\centerline{\graphicspath{ {Manifold_35/cancelling5_b-/} }\includegraphics[width=14cm, height=12cm]{all_x-y_x-z}}

\centerline{\graphicspath{ {Manifold_35/cancelling5_b-/} }\includegraphics[width=12cm, height=10cm]{all_y-z_x-y-z}}

We move on to cancelling the 1-handle $A-,A'-$ using the yellow 2-handle in the x-y plane. This cancellation affects the x-y plane
and the diagram corresponding to the twelve 2-handles that did not all lie in a single plane.

\centerline{\graphicspath{ {Manifold_35/cancelling6_a-/} }\includegraphics[width=10cm, height=5cm]{2cycle_x-y_int}}

\centerline{\graphicspath{ {Manifold_35/cancelling6_a-/} }\includegraphics[width=10cm, height=5cm]{2cycle_x-y-z}}

The intersection points in the y-z plane change as follows.

\centerline{\graphicspath{ {Manifold_35/cancelling6_a-/} }\includegraphics[width=11cm, height=6cm]{2cycle_y-z_int}}

We can then cancel $D,D'$ using the brown 2-handle in the diagram corresponding to the twelve 2-handles that did not all lie in a single plane.
The result of this cancellation can be seen in the following picture.

\centerline{\graphicspath{ {Manifold_35/cancelling7_d/} }\includegraphics[width=10cm, height=5cm]{2cycle_x-y-z}}

The cancellation also affects the x-z plane, which changes in the following way.

\centerline{\graphicspath{ {Manifold_35/cancelling7_d/} }\includegraphics[width=10cm, height=5cm]{2cycle_x-z_int}}

The intersection points in the x-y plane change in the following way.

\centerline{\graphicspath{ {Manifold_35/cancelling7_d/} }\includegraphics[width=10cm, height=5cm]{2cycle_x-y_int}}

In the x-z plane we have the blue 2-handle that runs over $D-,D'-$ once, these then form a handle cancellation pair.
Carrying out this cancellation, the x-z plane changes to:

\centerline{\graphicspath{ {Manifold_35/cancelling8_d-/} }\includegraphics[width=10cm, height=5cm]{2cycle_x-z_int}}

The diagram corresponding to the twelve 2-handles that did not all lie in a single plane will also be affected by this cancellation
and change in the following way.

\centerline{\graphicspath{ {Manifold_35/cancelling8_d-/} }\includegraphics[width=10cm, height=5cm]{2cycle_x-y-z}}

The intersection points labelled \textbf{C-D-} and \textbf{C'-D'-} will disappear from the y-z plane, no new intersection points
will appear.

\centerline{\graphicspath{ {Manifold_35/cancelling8_d-/} }\includegraphics[width=10cm, height=5cm]{2cycle_y-z_int}}

As for the x-y plane, we have that the intersection points labelled \textbf{XZ\_D-D'-} will disappear but some new ones
corresponding to the yellow 2-handles in the x-z plane, running from $K-$ to $L-$ and from $K'-$ to $L'-$, will appear. We will also
see intersection points created by the two 2-handles running between $C-$ and $C'-$ in the diagram corresponding to the twelve 2-handles that did not all lie in a single 
plane.

\centerline{\graphicspath{ {Manifold_35/cancelling8_d-/} }\includegraphics[width=10cm, height=5cm]{2cycle_x-y_int}}

We can then cancel $C-,C'-$ with the light green 2-handle in the x-z plane that runs over it once. Note that when we carry out this cancellation a few
2-handles can be immediately deleted from our diagrams. Namely, the orange 2-handle in the x-z plane will slide to give a zero framed unknot, and so
will the light green and yellow 2-handles in the  diagram corresponding to the twelve 2-handles that did not all lie in a single 
plane. These then each cancel a 3-handle and hence can be deleted from our diagrams.

The following picture shows how the x-z plane changes after this cancellation has been carried out.

\centerline{\graphicspath{ {Manifold_35/cancelling9_c-/} }\includegraphics[width=10cm, height=5cm]{2cycle_x-z_int}}

The diagram corresponding to the twelve 2-handles that did not all lie in a single plane changes in the following way.

\centerline{\graphicspath{ {Manifold_35/cancelling9_c-/} }\includegraphics[width=9cm, height=4cm]{2cycle_x-y-z}}

The intersection points labelled \textbf{XZ\_{C-C'-}} and \textbf{C-C'-} in the x-y plane will disappear. The following
picture shows the coding of the new ones that come into place.

\centerline{\graphicspath{ {Manifold_35/cancelling9_c-/} }\includegraphics[width=10cm, height=5cm]{2cycle_x-y_int}}

We move on to cancelling $J,J'$ with the dashed black 2-handle that resides in the x-y and y-z planes. The x-y plane changes in the following
way.

\centerline{\graphicspath{ {Manifold_35/cancelling10_j/} }\includegraphics[width=11cm, height=6cm]{2cycle_x-y_int}}

The reader should be aware that when we cancel $J,J'$ with the dashed black 2-handle we obtain grey and black 2-handle components that loop
back into $G'$ and $H$. We can then push these through the pieces of 1-handles they loop back into and then do a handle slide to obtain
grey and black 2-handles that run over $I-,I'-$ once, which can be seen in the above picture. 

The y-z plane changes in the following way.

\centerline{\graphicspath{ {Manifold_35/cancelling10_j/} }\includegraphics[width=11cm, height=6cm]{2cycle_y-z_int}}

The reader should be aware that, as in the case of the x-y plane, we have carried out some handle slides as well. Namely, the cancellation
creates a red 2-handle that loops back into $K'$ and a green 2-handle that loops back into $L'$, we can then slide these into the positions shown.
The dashed black 2-handle that we used to cancel $J,J'$ intersected the x-z plane, hence the cancellation will add some new points of intersection
with the x-z plane. The following two diagrams show these new intersection points with the second one being a close up showing the
coding of the intersection points in the middle of the diagram.

\centerline{\graphicspath{ {Manifold_35/cancelling10_j/} }\includegraphics[width=11cm, height=6cm]{2cycle_x-z_int}}

\centerline{\graphicspath{ {Manifold_35/cancelling10_j/} }\includegraphics[width=6cm, height=3cm]{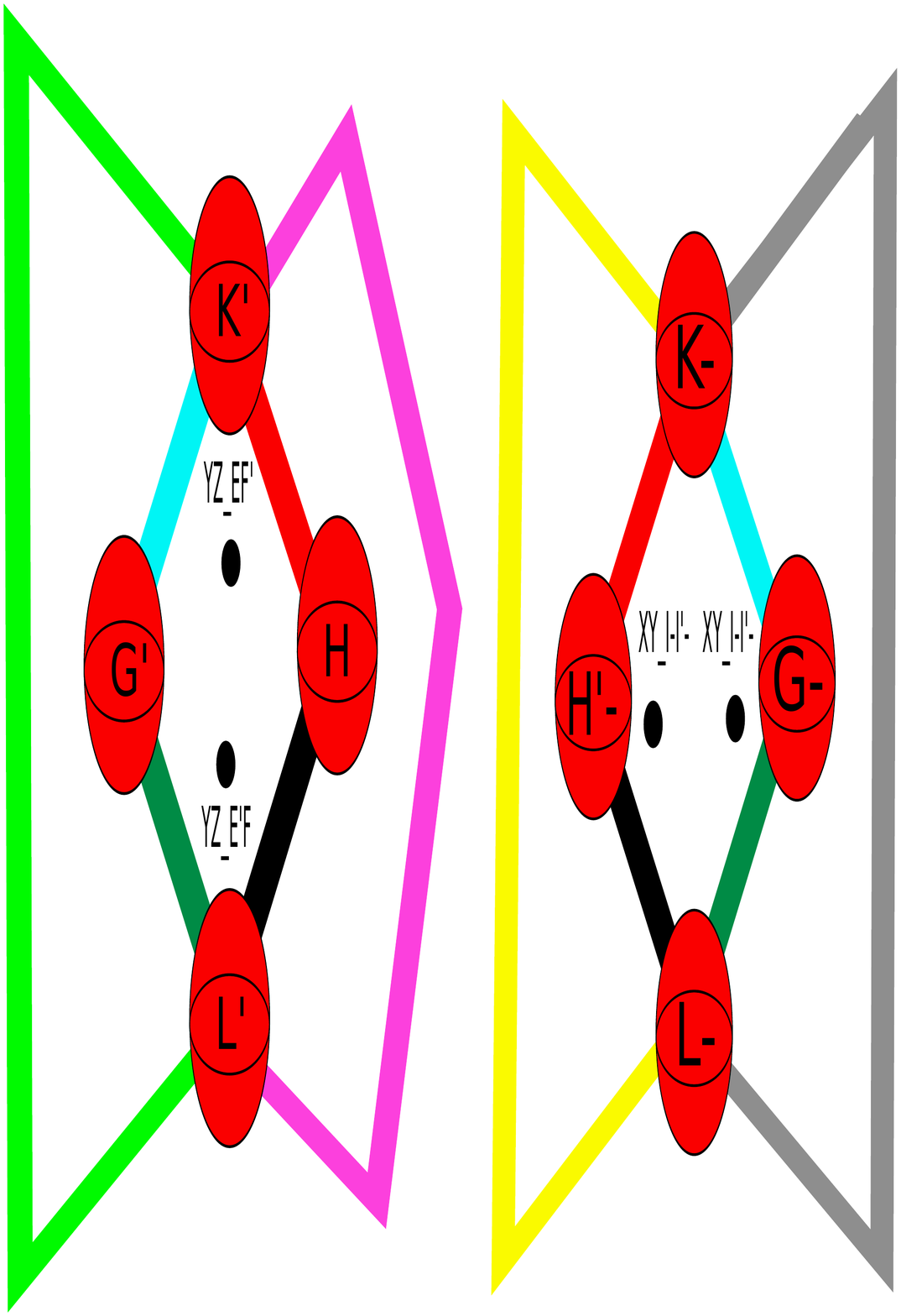}}

At this point we stop and collect together pictures of the four diagrams so far. This first picture shows the structure of the
x-y and x-z planes, and the second picture shows the structure of the y-z plane and the 2-handles remaining in the diagram corresponding to the twelve 2-handles that did 
not all lie in a single plane.

\centerline{\graphicspath{ {Manifold_35/cancelling10_j/} }\includegraphics[width=12cm, height=10cm]{all_x-y_x-z}}

\centerline{\graphicspath{ {Manifold_35/cancelling10_j/} }\includegraphics[width=10cm, height=8cm]{all_y-z_x-y-z}}

The next handle cancellation we carry out is to cancel $K,K'$ with the black dashed 2-handle that sits in the x-z and y-z planes.

The following pictures shows how the x-z and y-z planes change respectively.

\centerline{\graphicspath{ {Manifold_35/cancelling11_k/} }\includegraphics[width=11cm, height=5cm]{2cycle_x-z_int}}

\centerline{\graphicspath{ {Manifold_35/cancelling11_k/} }\includegraphics[width=11cm, height=5cm]{2cycle_y-z_int}}

The dashed black 2-handle that was used to cancel $K,K'$ did not intersect the x-y plane, hence on the level of intersection points the x-y plane does
not change.

We can also cancel $L,L'$ with the green 2-handle that runs over it once in the x-z plane. This cancellation will only affect the
x-z and y-z planes. As the 2-handle we are using to perform the cancellation does not intersect the x-y plane we find that the intersection
points of the x-y plane remain the same.

The following pictures show how these diagrams look like after the cancellation has taken place.

\centerline{\graphicspath{ {Manifold_35/cancelling12_l/} }\includegraphics[width=11cm, height=5cm]{2cycle_x-z_int}}

\centerline{\graphicspath{ {Manifold_35/cancelling12_l/} }\includegraphics[width=11cm, height=5cm]{2cycle_y-z_int}}

The reader should note that we have also carried out some handle slides. After performing the above cancellation we get an orange 2-handle component
in the y-z plane that loops back into $E'$, and a brown 2-handle component (in the y-z plane as well) that loops back into $F$. We can then
perform a handle slide on both these components to obtain orange and brown 2-handles, in the y-z plane, that run over $K-,K'-$ once.

We can use either of these 2-handles to perform a cancellation with $K-,K'-$, again this will only affect the x-z and y-z planes. The following
two pictures show how these planes change.

\centerline{\graphicspath{ {Manifold_35/cancelling13_k-/} }\includegraphics[width=10cm, height=5cm]{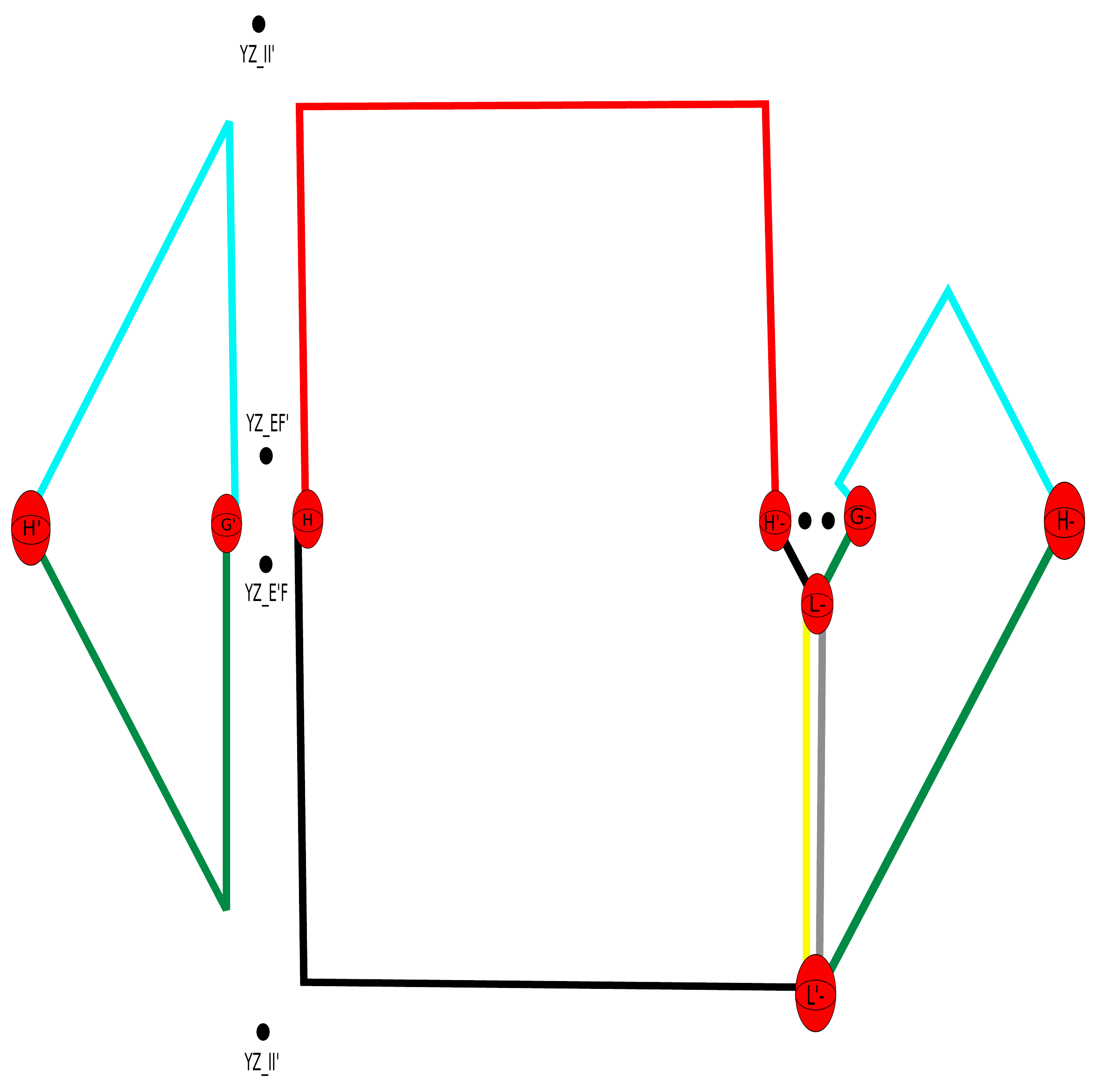}}

\centerline{\graphicspath{ {Manifold_35/cancelling13_k-/} }\includegraphics[width=10cm, height=5cm]{2cycle_y-z_int}}

We can also cancel $L-,L'-$ with the blue 2-handle in the y-z plane, this cancellation only affects the x-z and y-z planes.

The x-z plane changes to the following diagram.

\centerline{\graphicspath{ {Manifold_35/cancelling14_l-/} }\includegraphics[width=10cm, height=4cm]{2cycle_x-z_int}}

The y-z plane changes to the following diagram.

\centerline{\graphicspath{ {Manifold_35/cancelling14_l-/} }\includegraphics[width=11cm, height=5cm]{2cycle_y-z_int}}

This cancellation does not affect the points of intersection in the x-y plane as the 2-handle we used in the cancellation did not intersect the 
x-y plane.

The next cancellation we undertake is to cancel $I-,I'-$ with the black 2-handle in the x-y plane. This cancellation only affects
the 2-handles in the x-y and y-z plane.

The x-y plane changes to the following diagram:

\centerline{\graphicspath{ {Manifold_35/cancelling15_i-/} }\includegraphics[width=10cm, height=4cm]{2cycle_x-y_int}}

The y-z plane changes in the following way:

\centerline{\graphicspath{ {Manifold_35/cancelling15_i-/} }\includegraphics[width=11cm, height=6cm]{2cycle_y-z_int}}

The black 2-handle used to carry out this cancellation intersected the x-z plane, the code was \textbf{XY\_{I-I'-}}, therefore after the cancellation
this intersection point disappears. However, many new intersection points arise from 2-handle components in the y-z plane. The following picture
shows the coding of these new intersection points.

\centerline{\graphicspath{ {Manifold_35/cancelling15_i-/} }\includegraphics[width=11cm, height=6cm]{2cycle_x-z_int}}

We have carried out a further five handle cancellations. This is a good point to stop and take stock of how our four different diagrams,
showing the structure of all the 2-handles, look like.

The following picture shows the structure of the x-y and x-z planes respectively after all the above cancellations have been carried out.

\centerline{\graphicspath{ {Manifold_35/cancelling15_i-/} }\includegraphics[width=11cm, height=9cm]{all_x-y_x-z}}

The following picture shows the structure of the y-z plane and the structure of the diagram corresponding to the twelve 2-handles that did not all lie in a single plane
after all the above cancellations have been carried out.

\centerline{\graphicspath{ {Manifold_35/cancelling15_i-/} }\includegraphics[width=12cm, height=10cm]{all_y-z_x-y-z}}

We remind the reader that all throughout the above cancellations there was a 2-handle that we were not showing in our diagrams. Namely, the 2-handle corresponding
to the translation $e^{-1}heh^{-1}$. Recall that this 2-handle had four components running from $E$ to $H$, $H'-$ to $E-$, $E'$ to $H'$ and $H-$ to $E'-$.
So far, the cancellations we have carried out have all been within the three planes, the x-y, x-z, y-z planes, or the diagram showing the twelve 2-handles
not all lying in a single plane. The components of the 2-handle $e^{-1}heh^{-1}$ each pass between the y-z and x-y planes and so are not affected by any
of the cancellations that are carried within the x-y, x-z or y-z planes. As for cancellations done within in the diagram
corresponding to the twelve 2-handles that did not all lie in a single plane, two components, namely the ones running from $E$ to $H$
and $H'-$ to $E-$, are contained in the ``inside'' of the diagram, and it is easy to see that the cancellations we have carried out so far have
not in any way interfered with these two components. As for the two components running from $E'$ to $H'$ and $H-$ to $E'-$, these run on the 
``outside'' of the diagram, hence it is clear that the cancellations we have done so far have not interfered with these two components.

The next cancellation we are going to carry out is to cancel $I,I'$ using the red 2-handle that resides in the y-z plane. This cancellation
will affect the 2-handles in the y-z and x-y planes.

The y-z plane changes to the following diagram.

\centerline{\graphicspath{ {Manifold_35/cancelling16_i/} }\includegraphics[width=10cm, height=5cm]{2cycle_y-z_int}}

The x-y plane changes in the following way.

\centerline{\graphicspath{ {Manifold_35/cancelling16_i/} }\includegraphics[width=10cm, height=5cm]{2cycle_x-y_int}}

The 2-handle used to carry out the above cancellation intersected the x-z plane, the following diagram shows how the intersection points
in the x-z plane change.

\centerline{\graphicspath{ {Manifold_35/cancelling16_i/} }\includegraphics[width=10cm, height=5cm]{2cycle_x-z_int}}

We move on to cancelling $J-,J'-$ with the orange 2-handle in the x-y plane, this cancellation only affects the 2-handles in the x-y and y-z planes.

The x-y plane changes in the following way

\centerline{\graphicspath{ {Manifold_35/cancelling17_j-/} }\includegraphics[width=9cm, height=4cm]{2cycle_x-y_int}}

and the y-z plane changes to the following diagram.

\centerline{\graphicspath{ {Manifold_35/cancelling17_j-/} }\includegraphics[width=9cm, height=4cm]{2cycle_y-z_int}}

The orange 2-handle used in this cancellation intersected the x-z plane. The intersection points in the x-z plane changes to the following.

\centerline{\graphicspath{ {Manifold_35/cancelling17_j-/} }\includegraphics[width=10cm, height=5cm]{2cycle_x-z_int}}

So far, the last few cancellations we have carried out have not interfered with the 2-handles that are left in the 
diagram corresponding to the twelve 2-handles that did not all lie in a single plane. We want to perform an isotopy of the 2-handles in this diagram.

The first picture shows the original position the 2-handles were in, and the picture following it shows the final position after we carry out the isotopy.
It should be clear to the reader how the 2-handles move during this isotopy.

\centerline{\graphicspath{ {Manifold_35/cancelling18_isotopy/} }\includegraphics[width=9cm, height=8cm]{2cycle_x-y-z_2}}

The next cancellation we are going to undertake is to cancel $H,H'-$ using the red 2-handle in the x-z plane. This cancellation
affects the 2-handles in the x-z and x-y planes. It also affects the 2-handle corresponding to the translation $e^{-1}heh^{-1}$ that we have
not been drawing so far.

We start with the x-z plane, the following picture shows how the x-z plane changes after we have carried out this cancellation.

\centerline{\graphicspath{ {Manifold_35/cancelling19_h/} }\includegraphics[width=7cm, height=3cm]{2cycle_x-z_int}}

The x-y plane changes in the following way:

\centerline{\graphicspath{ {Manifold_35/cancelling19_h/} }\includegraphics[width=7cm, height=3.5cm]{2cycle_x-y_int}}

Recall that the 2-handle corresponding to the translation $e^{-1}heh^{-1}$ had in total four components, two of them in particular were such that
one ran from $E$ to $H$ and another from $E-$ to $H'-$. When we perform the above cancellation these two components come together, giving one component
running from $E$ to $E-$.

The following picture shows how this new component runs between $E$ and $E-$.

\centerline{\graphicspath{ {Manifold_35/cancelling19_h/} }\includegraphics[width=5cm, height=4cm]{2cycle_x-y-z_extra_handle}}

The cancellation carried out above causes the intersection points in the y-z plane to change. The following picture
shows the y-z plane with these new intersection points.

\centerline{\graphicspath{ {Manifold_35/cancelling19_h/} }\includegraphics[width=8cm, height=4cm]{2cycle_y-z_int}}

We can then cancel $H',H-$ with the green 2-handle in the x-y plane. This cancellation affects the 2-handles in the x-y plane and the x-z plane.

The x-y plane changes as follows.

\centerline{\graphicspath{ {Manifold_35/cancelling20_h-/} }\includegraphics[width=8cm, height=3cm]{2cycle_x-y_int}}

The x-z plane changes to the following.

\centerline{\graphicspath{ {Manifold_35/cancelling20_h-/} }\includegraphics[width=8cm, height=4cm]{2cycle_x-z_int}}

The translation $e^{-1}heh^{-1}$ that has been reduced to consisting of one component running from $E$ to $E-$, another running from 
$E'$ to $H'$, and another running from $E'-$ to $H-$. When we cancel $H',H-$ with the green 2-handle in the x-y plane the components
running from $E'$ to $H'$ and $E'-$ to $H-$ slide to form one component joining $E'$ to $E'-$.

The following picture shows how this new component runs between $E'$ and $E'-$.

\centerline{\graphicspath{ {Manifold_35/cancelling20_h-/} }\includegraphics[width=8cm, height=6cm]{2cycle_x-y-z_extra_handle}}

We now have two separate diagrams, the one coming from the x-y and x-z planes that involve 2-handles running over $G',G-$, and the 2-handles
left in the diagram corresponding to the twelve 2-handles that did not all lie in a single plane, and the 2-handles in the y-z plane, which run
over the 1-handles $E,E'-$, $E',E-$, $F,F'-$ and $F',F-$. It is easy to see that these two diagrams do not interact with each other in any way.
Furthermore, the 2-handles that run over $G',G-$ (the ones coming from the x-y and x-z planes) do so once. Hence we can use any one of them to form
a cancelling pair with $G',G-$. Carrying out this cancellation, all other 2-handles running over $G',G-$ slide to form zero framed unknot's, hence cancel
with a 3-handle and can be deleted from the diagram. Thus we are left with the 2-handles running over the 1-handles $E,E'-$, $E',E-$, $F,F'-$ and $F',F-$.

We can then cancel $F',F-$ with the black 2-handle in the y-z plane. This will cause the pink 2-handle in the y-z plane to slide and have two components, 
one looping back into $E$ and the other looping back into $E'-$. We can then slide one of these components through to obtain a zero framed unknot
that cancels a 3-handle. Hence this 2-handle can be deleted from the diagram. \\
We can also cancel $F,F'-$ using the dark green 2-handle in the y-z plane. This cancellation causes the yellow 2-handle in the y-z plane to slide into
a position where it has one component looping back into $E'$ and another looping back into $E-$. We can then perform a handle slide to obtain a yellow
coloured zero framed unknot. This then cancels with a 3-handle and can be deleted from the diagram. 

The two cancellations we have just carried out also affect the other 2-handles, they reside in the diagram corresponding to the twelve 2-handles that did not all lie in a 
single plane. It is straightforward to see how they change, the following diagram shows the position they slide into.

\centerline{\graphicspath{ {Manifold_35/cancelling22_f_and_f-/} }\includegraphics[width=9cm, height=4cm]{2cycle_x-y-z}}

We also have an extra 2-handle that we have not shown in the above, it is the 2-handle that corresponded to the translation $e^{-1}heh^{-1}$. It has
two components, one that runs from $E$ to $E-$ and another that runs from $E'$ to $E'-$. The following picture adds this 2-handle to the above
diagram.

\centerline{\graphicspath{ {Manifold_35/cancelling23_add_extra_handle/} }\includegraphics[width=7cm, height=4cm]{2cycle_x-y-z_2}}

We can now carry out some handle slides. We can slide the orange 2-handle along the turquoise 2-handle to obtain the following diagram.

\centerline{\graphicspath{ {Manifold_35/cancelling24_slides/} }\includegraphics[width=7cm, height=4cm]{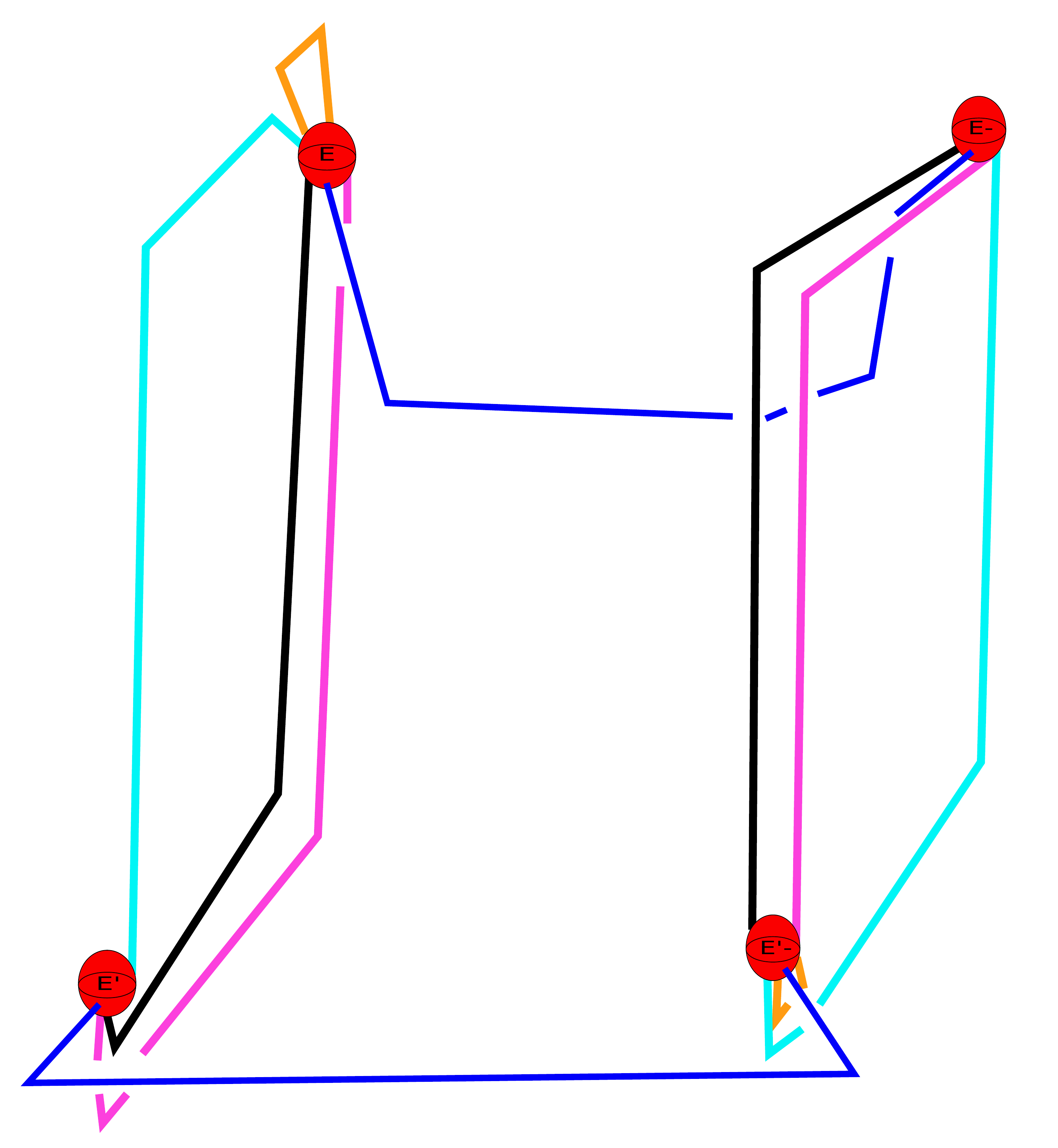}}

We then slide the pink 2-handle along the black 2-handle to obtain the following diagram.

\centerline{\graphicspath{ {Manifold_35/cancelling24_slides/} }\includegraphics[width=7cm, height=4cm]{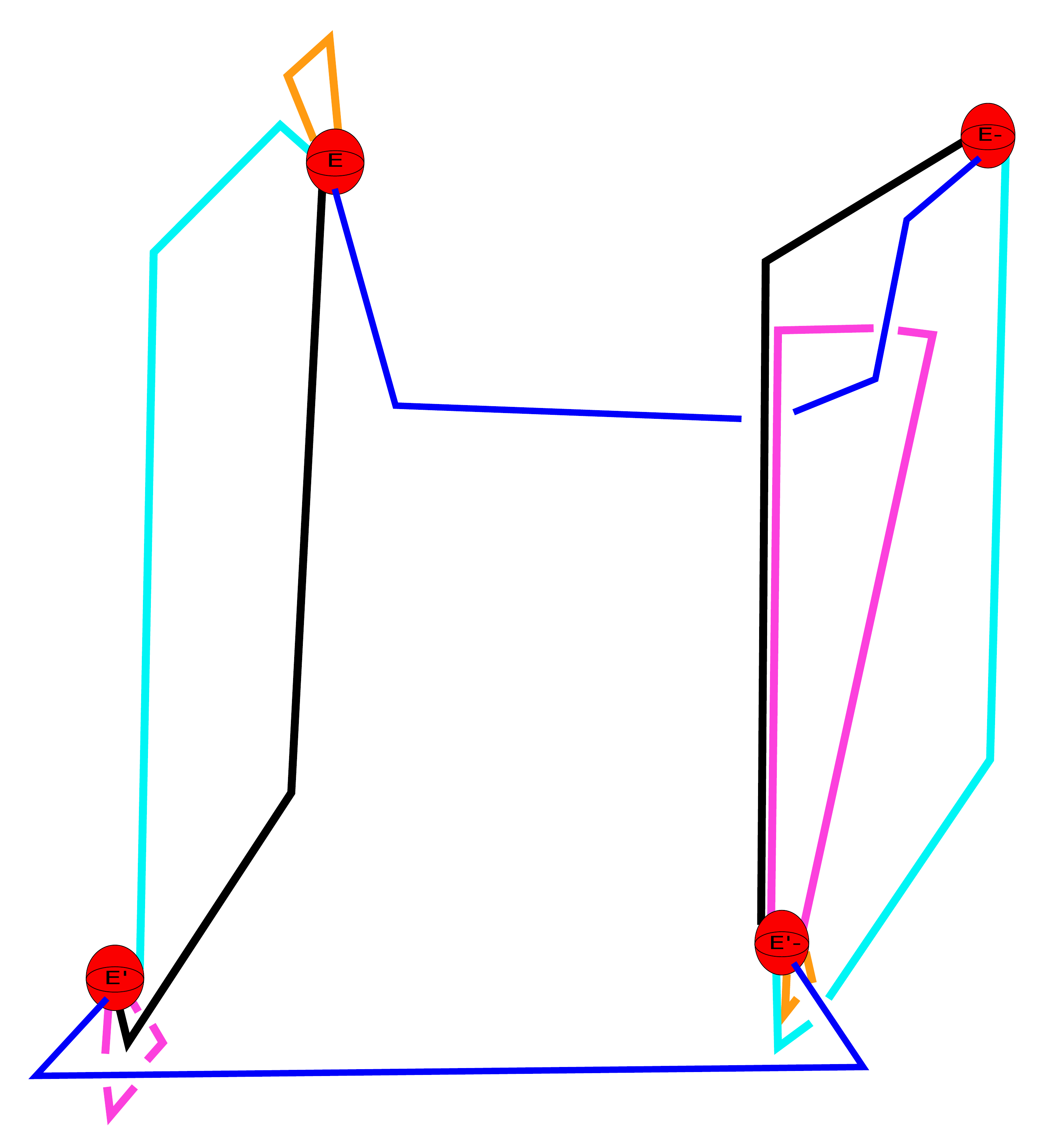}}

We can then slide the orange 2-handle to a zero framed unknot, which will then cancel a 3-handle. Therefore we can simply delete the orange 2-handle
from our diagram. This gives us the following diagram.

\centerline{\graphicspath{ {Manifold_35/cancelling24_slides/} }\includegraphics[width=7cm, height=4cm]{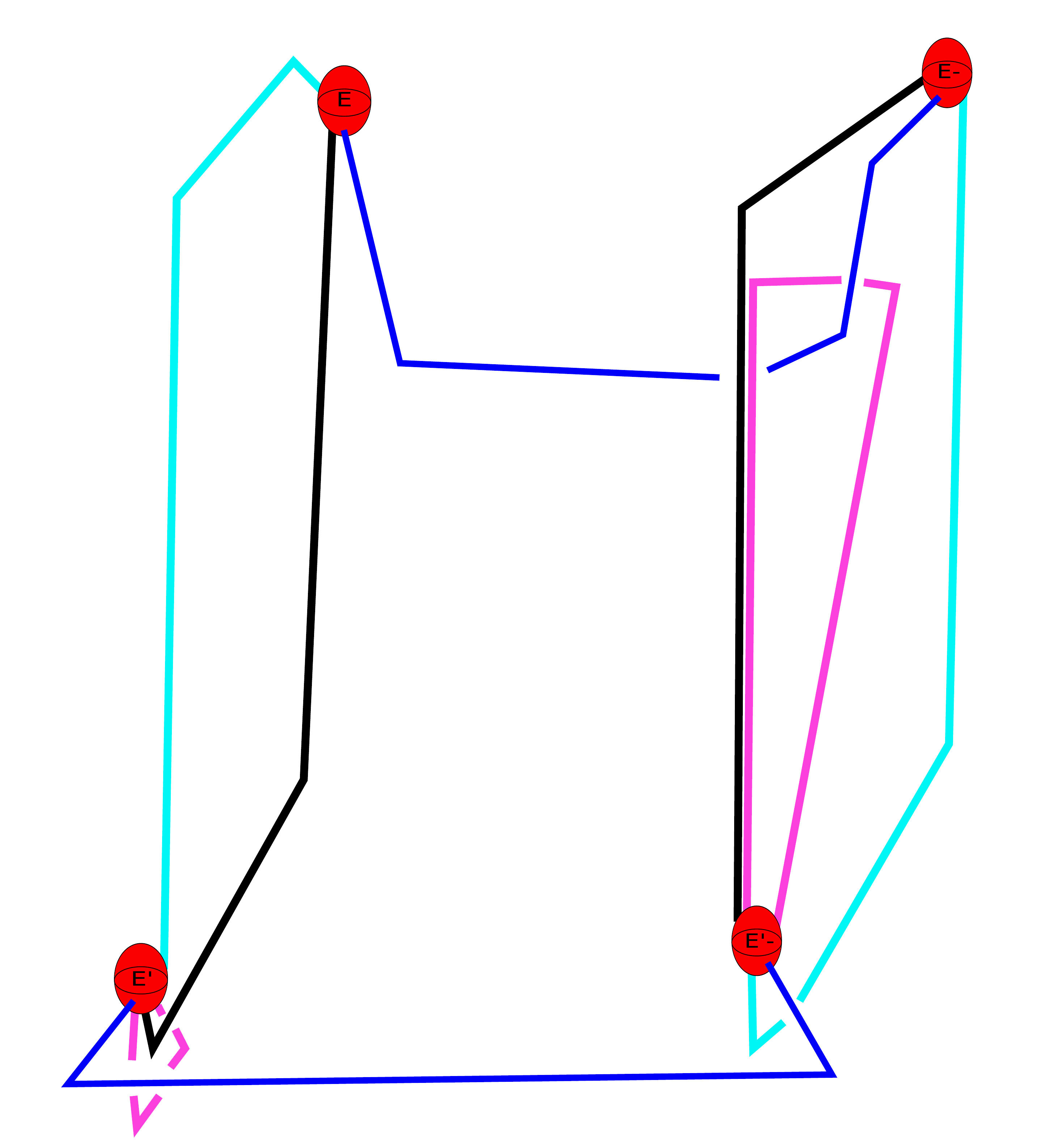}}

We can then slide the pink 2-handle into the following position.

\centerline{\graphicspath{ {Manifold_35/cancelling24_slides/} }\includegraphics[width=6cm, height=3cm]{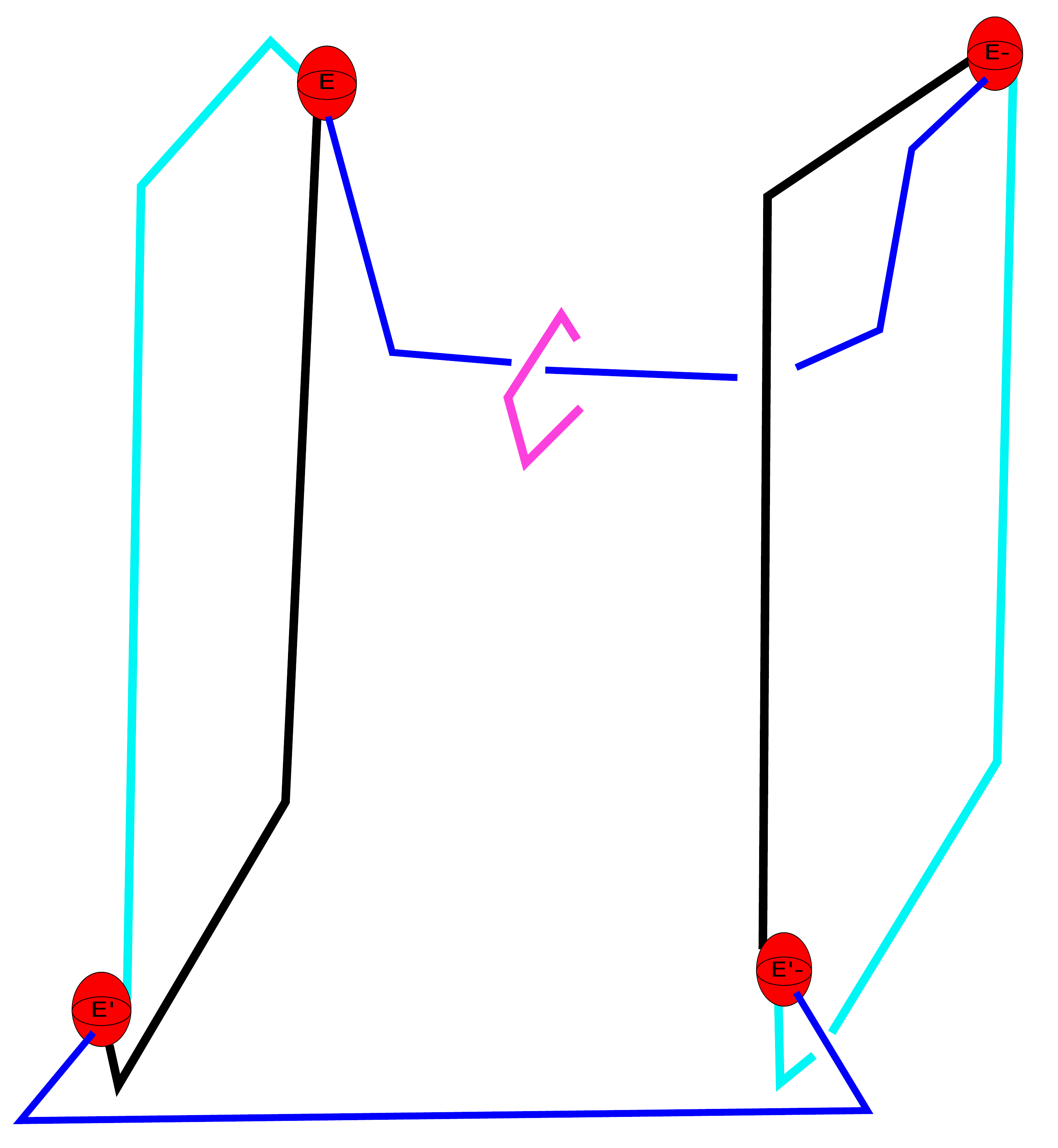}}

We can also slide the turquoise 2-handle along the black 2-handle to obtain the following diagram.

\centerline{\graphicspath{ {Manifold_35/cancelling25_e-/} }\includegraphics[width=6cm, height=3cm]{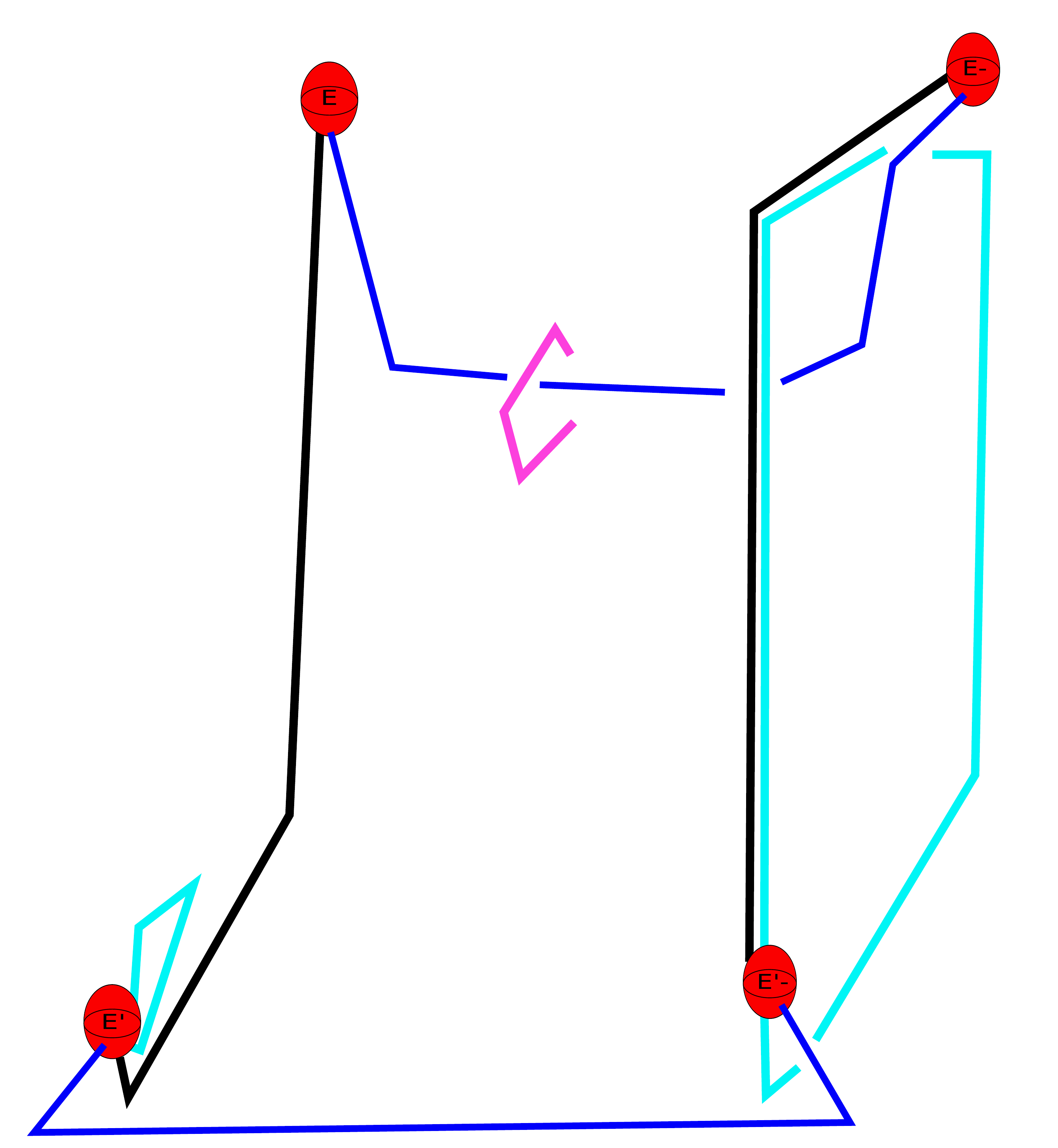}}

Then another handle slide produces:

\centerline{\graphicspath{ {Manifold_35/cancelling25_e-/} }\includegraphics[width=7cm, height=4cm]{2cycle_x-y-z_2}}

We can then cancel the 1-handle $E',E-$ with the black 2-handle producing the following diagram.

\centerline{\graphicspath{ {Manifold_35/cancelling25_e-/} }\includegraphics[width=5cm, height=3cm]{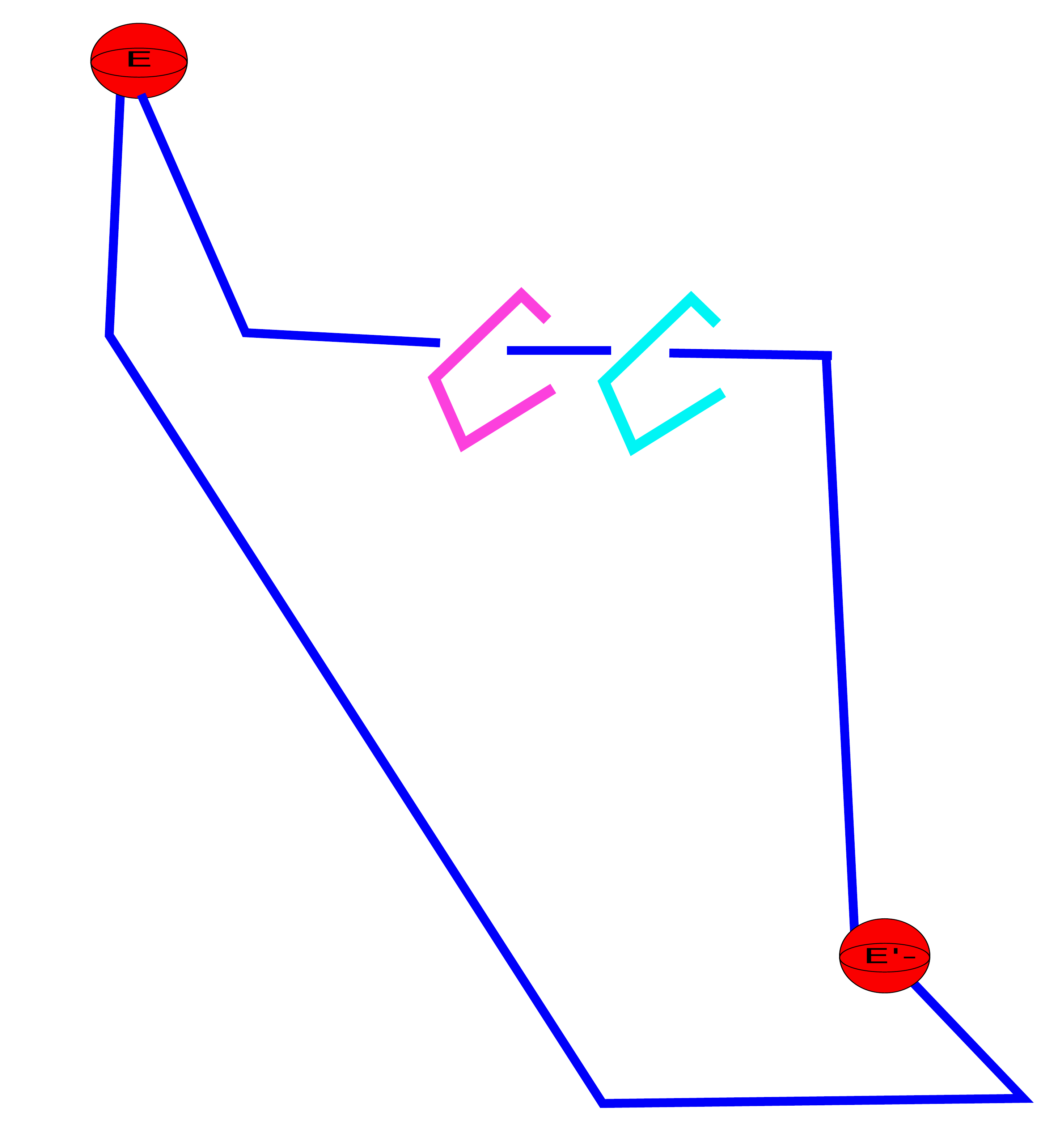}}

This is a Kirby diagram for the orientable manifold $\widetilde{M}$, which recall is the orientable double cover of $M$ (manifold 35). Furthermore, 
it is clear that $\pi_1(\widetilde{M}) = \langle x \hspace{0.2cm}  \vert \hspace{0.2cm}  x^2 = 1 \rangle \cong \Z_2$.  We want to take the double cover
of this manifold, which we denoted by $\widetilde{M}_2$.
Recall the procedure to do this, the one skeleton of $\widetilde{M}$ consists of $D^4 \cup (E,E'-)$ (the 0-handle union the 1-handle), which
is a copy of $S^1 \times D^3$. The one skeleton of $\widetilde{M}_2$ will also consist of a copy of $S^1 \times D^3$ double covering the one skeleton
of $\widetilde{M}$ in the usual way that $S^1 \times D^3$ double covers itself. Each of the remaining handles of $\widetilde{M}$ lift to two handles
of $\widetilde{M}_2$. This means that the blue 2-handle in the above diagram, that passes over $E,E'-$ twice, will lift to two 2-handles each passing
over the unique 1-handle in $\widetilde{M}_2$. The turquoise and pink 2-handles that loop around one component of the blue 2-handle in the above diagram, lift
to two copies of each looping around one lift of the blue 2-handle. The following diagram shows how the Kirby diagram of the double cover $\widetilde{M}_2$ looks
like.

\centerline{\graphicspath{ {Manifold_35/cancelling26_double_cover/} }\includegraphics[width=3cm, height=1cm]{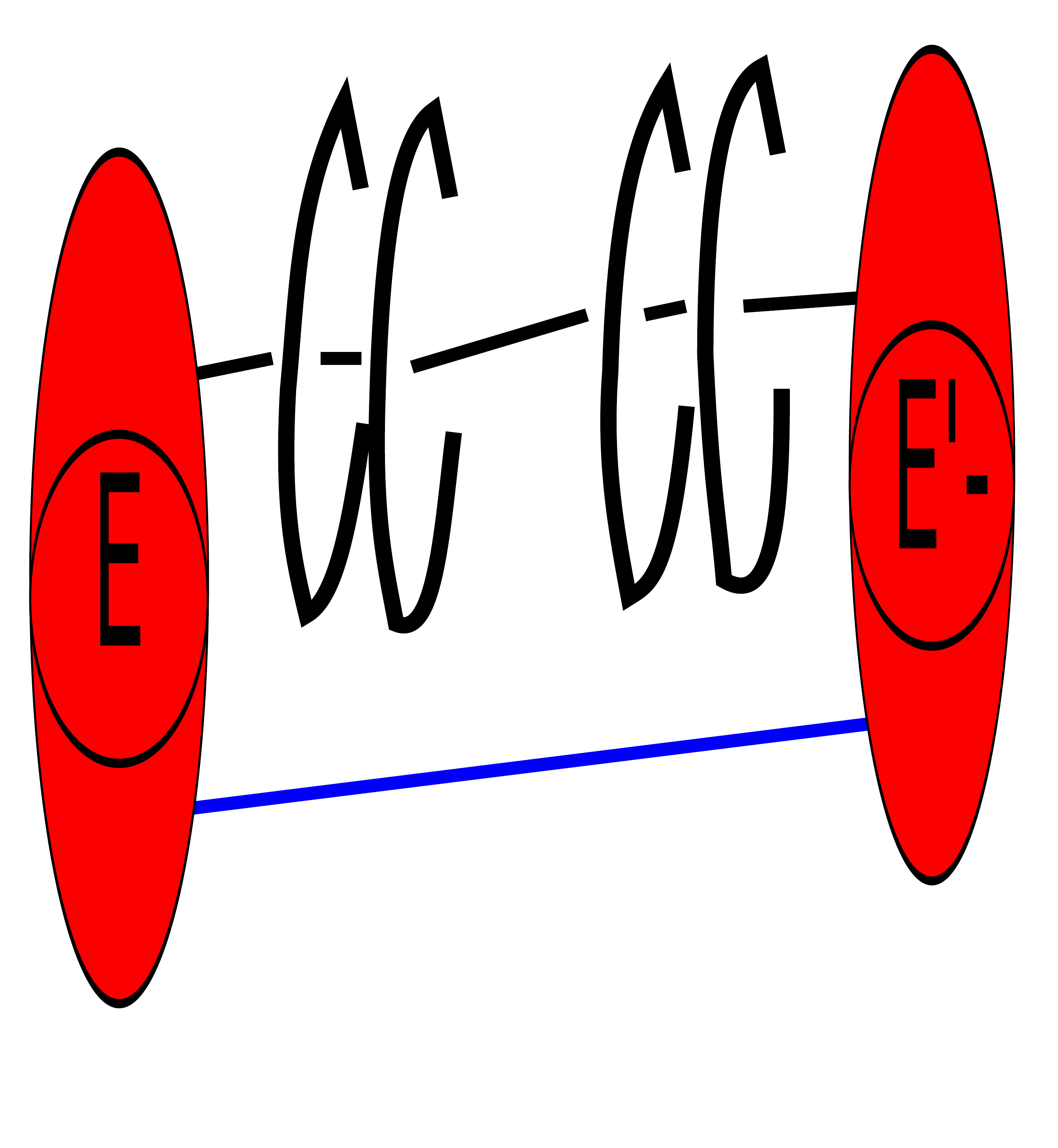}}

We can then cancel $E,E'-$ with the blue 2-handle, this will cause the diagram to change to.

\centerline{\graphicspath{ {Manifold_35/cancelling26_double_cover/} }\includegraphics[width=3cm, height=2cm]{2cycle_x-y-z_2}}

We can then slide three of the linked circles over the fourth one so that each one gives a zero framed unknot. These each cancel with
a 3-handle and we are left with the following diagram.

\centerline{\graphicspath{ {Manifold_35/cancelling26_double_cover/} }\includegraphics[width=3cm, height=2cm]{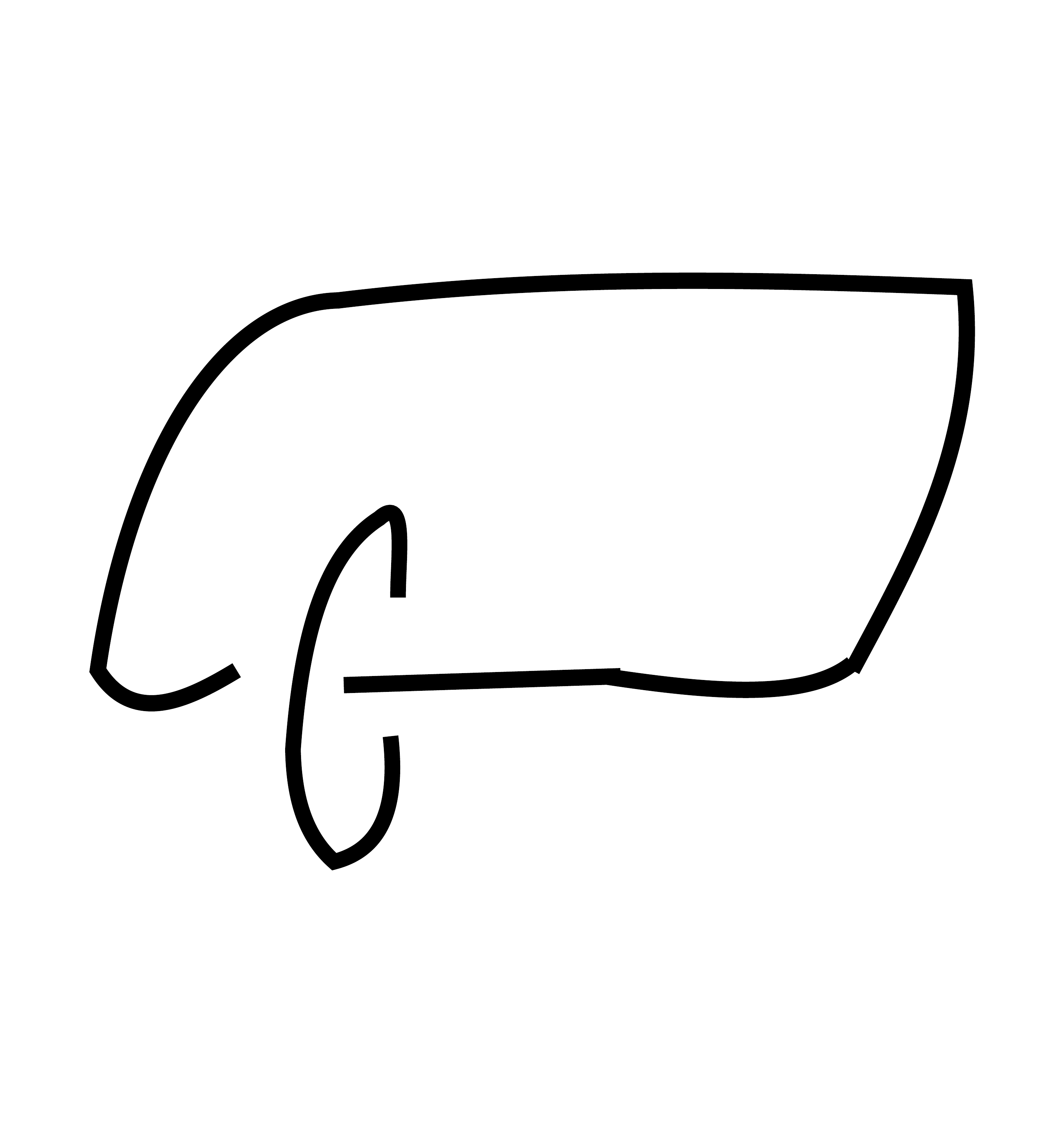}}

The diagram shows two zero framed linked 2-handles, which is precisely the Kirby diagram for $S^2 \times S^2$. Thus we can conclude that
the simply connected closed 4-manifold $\widetilde{M}_2$, which is a four fold cover of M (manifold 35), is diffeomorphic to $S^2 \times S^2$. 
We have thus proved the following theorem:

\begin{thm}\label{mainthm_2}
There exists a collection $L$ of linked 2-tori embedded in a standard $S^2 \times S^2$ such that the complement $(S^2 \times S^2) - L$
admits a finite volume hyperbolic geometry.
\end{thm}


\begin{thebibliography}{10}

%\bibitem{akbulut}
%Akbulut, S. \emph{4-manifolds 
%(2014)}, lecture notes available at http://www.math.msu.edu/~akbulut/.

%\bibitem{benedetti}
%Benedetti, R. \emph{Lectures on Hyperbolic
%Geometry}, Springer-Verlag, Berlin Heidelberg, 1992.


\bibitem{cox}
Coxeter, H.S.M \emph{Regular Complex Polytopes. 
 Second Edition.}, Cambridge University Press,
Cambridge, 1991.

\bibitem{gompf}
Gompf, R.E. and Stipsicz, A.I \emph{4-manifolds and Kirby
Calculus}, Graduate Studies in Mathematics, Providence, Rhode Island, 1999.


\bibitem{hantzsche}
Hantzsche, W. and Wendt, H. \emph{Dreidimensionale euklidische 
Raumformen}, Math. Ann. 110 (1935), 593-611.


%\bibitem{ivansic}
%Ivan$\check{s}$i$\acute{c}$, D. \emph{Hyperbolic structure on a 
%complement of tori in the 4-sphere}, Adv. Geom. 4, no. 1 (2004), 119–139.



\bibitem{johnson}
Johnson, D.L. \emph{Topics in the Theory of Group Presentations 
curvature}, London Mathematical Society Lecture Notes Series, 42, Cambridge, 1990.



\bibitem{kerckhoff}
Kerckhoff, S.P. and Storm, P.A. \emph{From the hyperbolic 24-cell
to the cuboctahedron}, Geometry \& Topology 
14 (2010) 1383–1477.

%\bibitem{kirby}
%Kirby, R. \emph{A Calculus for Framed Links 
%in $S^3$}, Inventiones math., 45 (1978), 35-56.


\bibitem{ratcliffe}
Ratcliffe, J.G. and Tschantz, S.T. \emph{The Volume Spectrum of
Hyperbolic 4-manifolds}, Experiment. Math.
Volume 9, Issue 1 (2000), 101-125.


\bibitem{sarat}
Saratchandran, H. \emph{Kirby diagrams and the Ratcliffe-Tschantz 
hyperbolic 4-manifolds}, ArXiv:math:GT/1503.06722.

\bibitem{sarat_2}
Saratchandran, H. \emph{A four dimensional hyperbolic link complement 
in a standard $S^4$}, ArXiv e-prints (2015).
 	 
%\bibitem{thurston_1}
%Thurston, W. \emph{The geometry and topology of 
%three-manifolds}, Princeton Univ. Math. Dept. Notes (1979), available at http://www.msri.org/communications/books/gt3m


%\bibitem{thurston_2}
%Thurston, W. \emph{Three-dimensional manifolds, Kleinian groups 
%and hyperbolic geometry}, Bull. Amer. Math. Soc. (N.S.)6, no. 3 (1982), 357–381.


%\bibitem{tschantz}
%Ivan$\check{s}$i$\acute{c}$, D. and Ratcliffe, J.G. and Tschantz, S.T. \emph{Complements of tori and Klein bottles in the 
%4-sphere that have hyperbolic structure}, Algebr. Geom. Topol. 5 (2005), 999–1026.



\bibitem{wolf}
Wolf, J.A. \emph{Spaces of constant 
curvature}, McGraw-Hill, United States of America, 1967.

  
 




  
  
\end{thebibliography}
\end{document}